\documentclass[final,5p,times]{elsarticle}
\usepackage{hyperref}
\usepackage{amsmath}
\usepackage{amssymb}
\usepackage{tabularx}
\usepackage{bm}
\usepackage[T1]{fontenc}
\usepackage{refcount}
\usepackage{upgreek}
\usepackage{csquotes}
\usepackage{mathrsfs}
\usepackage{soul}
\usepackage{comment}
\usepackage[usenames,dvipsnames,table]{xcolor}
\definecolor{jgreen}{HTML}{90bbbb}
\definecolor{myblue}{HTML}{7db3e8}
\definecolor{mybblue}{HTML}{93d6ec}
\definecolor{mydblue}{HTML}{6787e4}
\definecolor{mypurple}{HTML}{baa8f0}
\definecolor{mygreen}{HTML}{8cd9ad}
\definecolor{hyellow}{HTML}{ffe44d}
\definecolor{lred}{HTML}{8af48a}
\definecolor{Hred}{HTML}{f48a8a}
\definecolor{pastelpink}{HTML}{F8C8DC}
\definecolor{myorange}{HTML}{F7B787}
\definecolor{mypinkt}{HTML}{ec93bf}
\definecolor{raspberry}{HTML}{ee5db9}
\definecolor{retroorange}{HTML}{ffbf00}
\definecolor{retroyellow}{HTML}{ffd700}

\newcommand{\iu}{{i\mkern1mu}}

\begin{document}

\begin{frontmatter}
\ead{lidiajgomesdasilva.io, lidiajoana@pm.me}
\title{ \texttt{DiscoTEX 1.0}: \underline{Dis}continuous \underline{co}llocation and implicit-\underline{t}urned-\underline{ex}plicit (\texttt{IM\underline{TEX}}) integration \underline{symplectic}, symmetric numerical algorithms with high order jumps for differential equations II: time-integration extension to higher-orders of numerical convergence}

\author{Lidia J. Gomes Da Silva}

\affiliation{organization={School of Mathematical Sciences, Queen Mary University of London},
            city={London},
            postcode={E1 4NS}, 
            country={UK}}
\begin{abstract}
\texttt{DiscoTEX} is a highly accurate numerical algorithm for computing numerical weak-form solutions to distributionally sourced partial differential equations (PDE)s. The aim of this second paper, succeeding \cite{da2024discotex}, is to present its extension up to twelve orders. This will be demonstrated by computing numerical weak-form solutions to the distributionally sourced wave equation and comparing it to its exact solutions. The full details of the numerical scheme at higher orders will be presented. 
\end{abstract}

\end{frontmatter}

\section{Introduction}

\texttt{DiscoTEX} is a numerical algorithm which computes numerical weak-form solutions to distributionally sourced partial differential equations via the method of lines framework. Spatial discretisation is performed by correcting Lagrange integration formula via \underline{dis}continuous \underline{co}llocation methods and tim-integration is computed via implicit-\underline{t}urned-\underline{ex}plicit Hermite formulae corrected to the discontinuous case.\\ 

In the first paper \cite{da2024discotex}, we have shown how this algorithm could be used for the modelling of E(X)tremely and (E)xtreme Mass-Ratio-Inspirals (X)/(E)MRIs in comparison to previous time-domain numerical methods used. Highly accurate numerical solutions were obtained provided three numerical control factors are optimised via comprehensive numerical convergence tests. These user-specifiable optimisation factors are: 
\begin{enumerate}
    \item \textcolor{gray}{\texttt{[CTRL F\ref{controlfactor1}]} - Number of \texttt{N} nodes;}
    \item \textcolor{gray}{\texttt{[CTRL F\ref{controlfactor2}]} - Number of \texttt{J} jumps;} 
    \item \textcolor{gray}{\texttt{[CTRL F\ref{controlfactor3}]} - $\Delta \tau$, timestep size.} 
\end{enumerate}
Our results were further complemented by demonstrating \newline \texttt{DiscoTEX}'s superior performance against its closest relatives: the parent \texttt{DiscoIMP} - using the purely implicit Hermite scheme for numerical integration,  and \texttt{DiscoREX} - which uses the purely explicit \texttt{Runge-Kutta} scheme. Finally, we computed numerical solutions to the governing equations for a point particle on a circular geodesic in a Schwarzschild black hole with higher-order time-integration \texttt{IMTEX} Hermite schemes via \texttt{DiscoTEX}, see ref. \cite{da2024discotex} Table 8. From our results is clear, that with the current software implementation, increasing the algorithm's order did not change the accuracy within the same numerical setup, while taking significantly longer as we increased the order of numerical integration. \\ 

In this work, we provide the full details of the \texttt{DiscoTEX} (and relatives) algorithm for higher orders of numerical convergence via the computation and comparison of the numerical solutions to the distributionally-source wave equation which possesses exact solutions \cite{field2009discontinuous, field2010persistent, field2022discontinuous, field2023}. All the machinery necessary for the inclusion of higher-order time jumps emerging from a point particle prescribing a time-dependent trajectory is detailed. It is further observed the convergence rates are as expected for higher-order time-integration schemes and thus, the work may find application in problems with long evolutions that don't require a second interpolation step (computed by the \texttt{DiscoTEX} \enquote{interpolator}) such as XMRIs waveform modelling and black hole spectroscopy \cite{afshordi2023waveform}.\\

The paper is organised as follows. In Section \ref{discotex_nutshell} we will describe the discontinuous time-integration method at order-12 and compare it to other algorithms as used in \cite{da2024discotex, da2023hyperboloidal, phdthesis-lidia}. In Section \ref{higher_at_discotex} we will compute the numerical weak-form solutions to distributionally sourced wave equations with \texttt{DiscoTEX} at higher orders and compare all algorithms. Finally, we conclude in Section \ref{finale} discussing potential applications for higher order \texttt{DiscoTEX}'s schemes. 

\section{The \texttt{DiscoTEX} numerical recipe in a nutshell}\label{discotex_nutshell}

In this work we will demonstrate how to compute numerical weak-form solutions \cite{da2024discotex, field2009discontinuous, field2010persistent, field2022discontinuous, field2023} to distributionally sourced PDEs of Type I described in the hyperboloidal coordinate chart $(t,x) \rightarrow (t(\tau,\sigma), x(\sigma))$ as, 
\begin{align}
    & \textrm{Type I:} \ \ \square \Psi[\tau,\sigma] = F(\tau) \delta'(\sigma - \xi_{p}) + G(\tau)\delta(\sigma - \xi_{p}), 
    \label{distr_wavequation}
\end{align}
where  $\mathcal{\theta} = \gamma^{2}(t - x \, \dot{\xi}_{p} - |x - x_{p}|)$, $G(\tau) = \cos(\tau)= -\iu F(\tau)$ and $\rm{sgn}$ refers to the signum function with the $(t,x)$ coordinates mapped into the new chart. As showed by \cite{field2010persistent} (albeit in the $(t,x)$ coordinate chart) this admits the following exact solution: 

\begin{align}
    \label{ch2_distr_wavequation_delta1_plus_delta0}
    & \textrm{Type I:} \ \ \square \Psi(t,x) = G(t) \delta(x - r_{p}) + F(t)\delta'(x - r_{p}), \\
    & \Psi(t,x) = -\frac{1}{2} \sin \mathcal{\theta} + \frac{1}{2} \iu \gamma^{2} [\dot{r}_{p} + \text{sgn}(x-r_{p})] \cos \mathcal{\theta}. \label{ch2_distr_wavequation_delta1_plus_delta0_sol}
\end{align}

These equations admit numerical weak-form solutions via the ansatz: 
\begin{equation} 
    \Psi(\tau,\sigma)=  \Psi^{+}(\tau,\sigma) \Theta[\sigma-\xi_{p}(\tau)] + \Psi^{-}(t,r)\Theta[\xi_{p}(\tau)-\sigma], 
    \label{weakformsolution_wave}
\end{equation}
where $\Theta(z)$ is the Heaviside function defined as
\begin{equation}\Theta(z)= \left\{ 
\begin{array}{ccc}
1 &\text{for}& z>0, \\
\frac{1}{2}&\text{for}& z=0, \\
0 &\text{for}& z<0.  
\end{array}
\right. 
\label{ch2_HeavisideStepFunction}
\end{equation}
A time-dependent jump is,
\begin{align}
     &[[\Psi]](\tau) =\Psi^{+}(\tau,\xi_{p}(\tau)) - \Psi^{-}(\tau,\xi_{p}(\tau)).
     \label{ch2_timedependent_jump}
\end{align}

\subsection{Method-of-lines framework}\label{sec_mol}

The first step to solve equation \ref{weakformsolution_wave}, numerically, is to start by reducing the problem to a coupled system of first-order ordinary differential equations in time. The numerical problem reduces to a \texttt{1+1D} problem illustrated by 

\begin{equation}
    \partial_{\tau} \textbf{U} = \textbf{L} \cdot \textbf{U} + \mathcal{S}, \ \ \ \ \textbf{U} = \begin{pmatrix}\Psi\\\Pi
    \end{pmatrix},
    \label{ch3_general_pdes}
\end{equation}
where $\textbf{L}$ is a spatial differential operator and $\Psi(\tau,\sigma)$, $\Pi(\tau,\sigma)$ are the field variables and their partial derivatives with respect to time. \texttt{DiscoTEX} tackles this system by correcting Lagrange and Hermite interpolation formulae by building discontinuous discretisations in both the spatial and temporal directions. In the following sections, we briefly detail this and refer the reader to \cite{da2024discotex} where they have been comprehensively explained. 

\subsection{Spatial discretisation via discontinuous collocation methods}\label{sec2_disco_space}

$\textbf{U}(\tau,\sigma)$, can be discretised in space such that, effectively the numerical system to be solved is effectively time-dependent, as $\textbf{U}(\tau,\sigma) \rightarrow \textbf{U}(\tau, \sigma_{i}):= \textbf{U}_{i}(\tau)  = \textbf{U}(\tau)$ where $a\leq \sigma_{i} \leq b $ with the collocation nodes ranging from $0<i<N$. 

Essentially, we admit the weak-form solution given in equation \eqref{weakformsolution_wave} and work with the generic interpolating piecewise polynomial 
\begin{equation}
    p(\sigma) =  \sum^{N}_{j=0} \bigg[ \Psi_{j} + \Delta_{\Psi}(\sigma_{j} - \xi_{p}; \sigma- \xi_{p}) \bigg] \pi_{j}(\sigma), 
\end{equation}
where $\pi(\sigma)$ is the Lagrange basis polynomial (LBP) given as
\begin{equation}
 \pi_{j}(\sigma) = \prod^{N}_{\substack{k=0,\\ k\neq j}} \frac{\sigma - \sigma_{k}}{\sigma_{j} - \sigma_{k}}, 
  \label{ch3_sb2_LBP}
\end{equation}
and the $\Delta_{\Psi}$ function is given by
\begin{small}
\begin{align}
    &\Delta_{\Psi}(\sigma_{j} - \xi_{p}; \sigma - \sigma_{p}) = \nonumber \\
    &= \bigg[\Theta(\sigma_{i} - \xi_{p}) - \Theta(\sigma_{j} - \xi_{p})\bigg]g(\sigma_{j} - \xi_{p})  \; \rm{when} \; \sigma = \sigma_{i}. \ \ \ \ \ 
    \label{ch3_sb2_delta}    
\end{align}
\end{small}
\noindent The LBP is built using suitable collocation nodes as it was explained in \cite{da2024discotex, da2023hyperboloidal}. In the work that follows we will use Chebyshev pseudospectral collocation nodes as defined in Section 3 of \cite{da2024discotex}. The function $g(\sigma_{j} - \xi_{p})$ is defined as 
\begin{equation}
    g(\sigma_{j} - \xi_{p}) = \sum^{M}_{m=0} \frac{J_{m}}{m!}(\sigma_{j} - \xi_{p})^{m}
    \label{ch3_sb2_gVector_expSum}. 
\end{equation}
Our master field variables $\Psi$ and $\Pi$, as defined in equation \eqref{ch3_general_pdes},is then effectively approximated as 
\begin{equation}
    \Psi(\tau,\sigma) \approx \sum^{N}_{j=0} \bigg[ \Psi_{j}(\tau) + \Delta_{\Psi}(\sigma_{j} - \xi_{p}(\tau); \sigma - \xi_{p}(\tau))  \bigg] \pi_{j}(\sigma).
    \label{ch3_sb2_WFsolution}
\end{equation}
To be precise, all the differential operators in equation \eqref{distr_wavequation} and further specified in equations (\eqref{ch34_sb2_diffOperators_gamma}-\eqref{ch34_sb2_diffOperators_iota}), are discretised as, 
\begin{eqnarray}
  \partial_{\sigma}^{n}(\tau,\sigma)\Psi|_{\sigma = \sigma_{i}} = p^{(n)}(\sigma) = \sum^{N}_{j=0}  D^{(n)}_{ij}  \Psi_{j} + s_{i}^{(n)}(\tau), \ \ \ \ 
\label{cha2_spatial_disc_generic}
\end{eqnarray}
where $  s^{(n)}_{i}(\tau)$ is given as 
\begin{equation}
    s^{(n)}_{i}(\tau) = \sum^{N}_{j=0}D^{(n)}_{ij} \Delta_{\Psi}\big(\sigma_{j} - \xi_{p}(\tau); x_{i} - \xi_{p}(\tau)\big)  
    \label{ch3_sb2_SpaceSpource}
\end{equation}
and the user-specifiable high-order jumps in equation \eqref{ch3_sb2_gVector_expSum} obtained through the computation of the higher-order recurrence relation given as 
\begin{align}
    &J_{m+2}(\tau) =  -\bar{\gamma}^{2} \bigg[ \sum^{m}_{k=0} {m \choose k} \bigg(\upvarepsilon^{(k)}(\xi_{p}) \dot{J}_{m+1-k} + \iota^{(k)}(\xi_{p}) J_{m+1-k}  \nonumber \\
    &- V^{(k)}(\xi_{p}) J_{m-k} + \upvarrho^{(k)} (\xi_{p}) (\dot{J}_{m-k} - \dot{\xi}_{p} J_{m+1-k}) \nonumber \\ 
    &+ \Gamma^{(k)}(\ddot{J}_{m-k} - 2 \dot{\xi}_{p} \dot{J}_{m+1-k}- \ddot{\xi_{p}}J_{m+1-k}  ) \bigg)   \nonumber \\
    &+  \sum^{m}_{k=1} {m \choose k} \bigg(\dot{\xi}_{p}^{2} \Gamma^{(k)}(\xi_{p}) + \dot{\xi}_{p} \upvarepsilon^{(k)}(\xi_{p}) + \upchi^{(k)}(\xi_{p}) \bigg) J_{m+2-k} \bigg].
\label{rec_hyperboloidal}
\end{align}
where here for simplicity the time dependence on the RHS of the jumps has been suppressed and we have, $\bar{\gamma}^{-2} = \big(  \dot{\xi}_{p}^{2} \Gamma(\xi_{p}) -  \dot{\xi}_{p} \upvarepsilon(\xi_{p}) - \upchi(\xi_{p})\big)$. This machinery introduces two user-specifiable numerical optimisation control-factors that will need to be subjected to a comprehensive numerical convergence study \cite{da2024discotex, da2023hyperboloidal, phdthesis-lidia}, 
\begin{enumerate}
   \item \texttt{[CTRL F1]} - Number of nodes, \texttt{N}; \label{controlfactor1}
\item \texttt{[CTRL F2]} - Number of jumps, \texttt{J}.\label{controlfactor2}
\end{enumerate}

\subsection{Discontinuous time-integration via Hermite interpolation at higher-orders}\label{sec2_disco_time}
\begin{figure*}
    \centering
    \includegraphics[width=16cm]{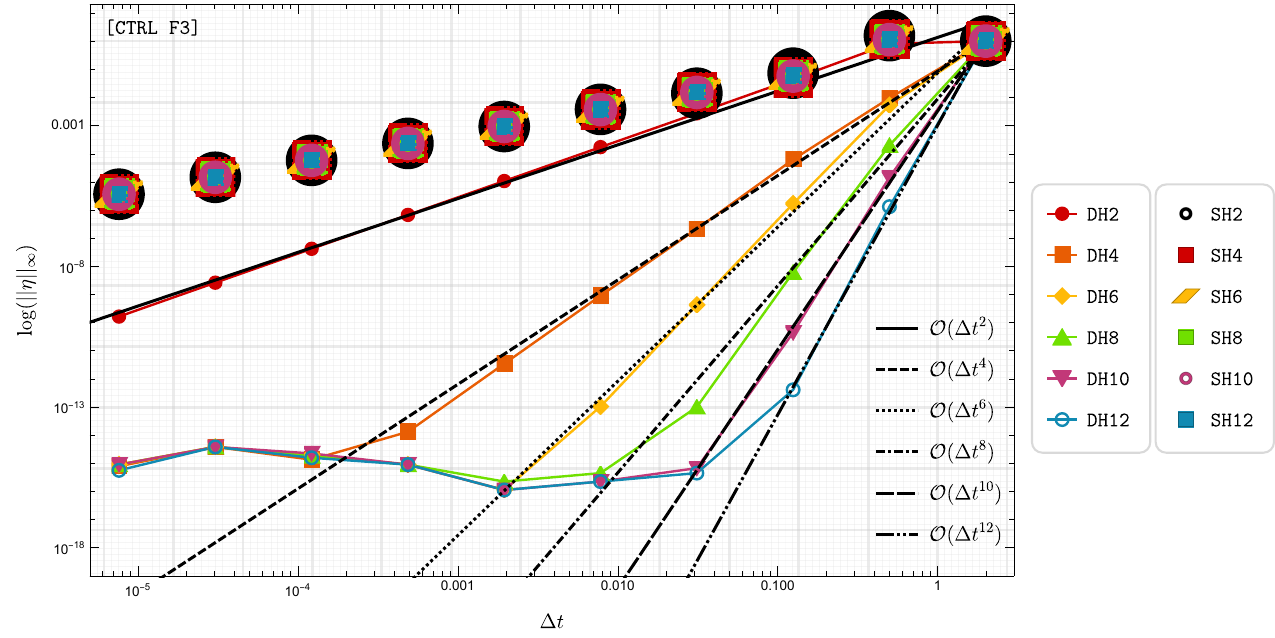}
    \caption{Numerical error associated with the numerical evaluation of the integral in equation \eqref{ch3_legendreQP_time} with the numerical scheme of equations (\eqref{app3_disco_time_h2} - \eqref{app3_disco_time_Jh10}) corrected by incorporating the discontinuous behaviour through the discontinuous time-integration rules given through equations (\eqref{app3_disco_time_Jh2} - \eqref{app3_disco_time_Jh12}). For numerical integration schemes of \texttt{2nd} to \texttt{12th order} numerical convergence in the same order is observed as given by the lines \texttt{DH2-DH12}. Inaccurate results are recorded for a smooth integrator as given by line \texttt{SH2-SH12}. The computational wall-clock times for the simulations are given in Table \ref{tab_discontinuousTime-allorders}. }
    \label{ch2_Disco_Time_AllOrders}
\end{figure*}

We now extend the procedure described in \cite{da2024discotex} to order-12 and  start by considering the first-order differential equation, 
\begin{equation}
    \frac{dy}{dt} = f(t, y(t)), 
    \label{ch2_disco_time_1st_ode}
\end{equation}
on a small time interval $[t_\nu, t_{\nu+1}]$. Applying the fundamental theorem of calculus we have, 
\begin{equation}
    y(t_{\nu+1}) - y(t_{\nu}) = \int_{t_{\nu}}^{t_{\nu+1}} f(t, y(t)) dt. 
    \label{ch2_disco_time_fundamental_theo_calculus}
\end{equation}
As demonstrated in \cite{da2024discotex} the discontinuous behaviour is incorporated by constructing the interpolant as a piecewise polynomial, 
\begin{equation}
    f(t, y(t)) \approx p(t) = p_{+}(t) \Theta(t - t_{\times}) +  p_{-}(t) \Theta(t_{\times} - t),  
    \label{ch2_disco_time_NS_interpo}
\end{equation}
where $t_{\times}$ is the point where the function is discontinuous such that $t_{\times} \in [t_{\nu}, t_{\nu+1}]$ and $f$ is approximated as $p_{+}$ in $[t_{\times}, t_{\nu + 1}]$ and $p_{-}$ in $[t_{\nu}, t_{\times}]$. Explicitly we have, 
\begin{align}
    \label{ch2_disco_time_system_right}
    &p_{+}(t) = a_{0} + a_{1} t  + a_{2} t^{2}  + a_{3} t^{3}  + a_{4} t^{4} + a_{5} t^{5} \nonumber \\
    &+ a_{6} t^{6}  + a_{7} t^{7}  + a_{8} t^{8} + a_{9} t^{9} + a_{10} t^{10}  + a_{11} t^{11},\\
    &p_{-}(t) = b_{0} + b_{1} t  + b_{2} t^{2}  + b_{3} t^{3}  + b_{4} t^{4} + b_{5} t^{5} \nonumber \\
    &+ b_{6} t^{6}  + b_{7} t^{7}  + b_{8} t^{8} + b_{9} t^{9} + b_{10} t^{10}  + b_{11} t^{11}.\\
    \label{ch2_disco_time_system_left}
\end{align}
Even though we have 24 unknown coefficients, the jump conditions are known and thus the following collocation conditions are imposed: 

\begin{align}
    &p_{-}(t_\nu) = f_{\nu},  p_{+}(t_{\nu+1}) = f_{\nu+1}, \\
    &p_{-}'(t_\nu) = df_{\nu},  \hspace{0.2cm} p_{+}'(t_{\nu+1}) = df_{\nu+1}, \\
    &p_{-}''(t_\nu) = ddf_{\nu},  \hspace{0.2cm} p_{+}''(t_{\nu+1}) = ddf_{\nu+1}, \\
    &p_{-}'''(t_\nu) = dddf_{\nu},  \hspace{0.2cm} p_{+}'''(t_{\nu+1}) = dddf_{\nu+1}, \\
    &p_{-}^{(4)}(t_\nu) = d4f_{\nu},  \hspace{0.2cm} p_{+}^{(4)}(t_{\nu+1}) = d4f_{\nu+1}, \\
    &p_{-}^{(5)}(t_\nu) = d5f_{\nu},  \hspace{0.2cm} p_{+}^{(5)}(t_{\nu+1}) = d5f_{\nu+1}, \\
    &p_{-}^{(6)}(t_\nu) = d6f_{\nu},  \hspace{0.2cm} p_{+}^{(6)}(t_{\nu+1}) = d6f_{\nu+1}, \\
    &p_{-}^{(7)}(t_\nu) = d7f_{\nu},  \hspace{0.2cm} p_{+}^{(7)}(t_{\nu+1}) = d7f_{\nu+1}, \\
    &p_{-}^{(8)}(t_\nu) = d8f_{\nu},  \hspace{0.2cm} p_{+}^{(8)}(t_{\nu+1}) = d8f_{\nu+1}, \\
    &p_{-}^{(9)}(t_\nu) = d9f_{\nu},  \hspace{0.2cm} p_{+}^{(9)}(t_{\nu+1}) = d9f_{\nu+1}, \\
    &p_{-}^{(10)}(t_\nu) = d10f_{\nu},  \hspace{0.2cm} p_{+}^{(10)}(t_{\nu+1}) = d10f_{\nu+1}, \\
    &p_{-}^{(11)}(t_\nu) = d11f_{\nu},  \hspace{0.2cm} p_{+}^{(11)}(t_{\nu+1}) = d11f_{\nu+1}, \\
    &p_{+}(t_\times)-  p_{-}(t_\times) = \textbf{J}_{0}, \\
    &p_{+}'(t_\times)-  p_{-}'(t_\times) = \textbf{J}_{1}, \\
    &p_{+}''(t_\times)-  p_{-}''(t_\times) = \textbf{J}_{2}, \\
    &p_{+}'''(t_\times)-  p_{-}'''(t_\times) = \textbf{J}_{3}, \\
    &p_{+}^{(4)}(t_\times)-  p_{-}^{(4)}(t_\times) = \textbf{J}_{4}, \\
    &p_{+}^{(5)}(t_\times)-  p_{-}^{(5)}(t_\times) = \textbf{J}_{5}, \\
    &p_{+}^{(6)}(t_\times)-  p_{-}^{(6)}(t_\times) = \textbf{J}_{6}, \\
    &p_{+}^{(7)}(t_\times)-  p_{-}^{(7)}(t_\times) = \textbf{J}_{7}, \\
    &p_{+}^{(8)}(t_\times)-  p_{-}^{(8)}(t_\times) = \textbf{J}_{8}, \\
    &p_{+}^{(9)}(t_\times)-  p_{-}^{(9)}(t_\times) = \textbf{J}_{9}, \\
    &p_{+}^{(10)}(t_\times)-  p_{-}^{(10)}(t_\times) = \textbf{J}_{10}, \\
    &p_{+}^{(11)}(t_\times)-  p_{-}^{(11)}(t_\times) = \textbf{J}_{11}. 
    \label{ch2_disco_time_24_collocation_condos}
\end{align}

With these 24 conditions, we can now solve for the 24 polynomial coefficients highlighted in equations (\eqref{ch2_disco_time_system_right}, \eqref{ch2_disco_time_system_left}) as a linear system of algebraic equations. Integrating both of the piecewise polynomials yields
\begin{align}
        &f(t)_{\texttt{DH12}} =\frac{\Delta t }{2} \bigg(f(t_{n}) + f(t_{n+1} )  \bigg) + \frac{5\Delta t^{2} }{44}  \bigg(\dot{f}(t_{n}) - \dot{f}(t_{n+1} )  \bigg) \nonumber \\
        &+ \frac{\Delta t^{3} }{66}  \bigg(\ddot{f}(t_{n}) + \ddot{f}(t_{n+1} )  \bigg) + 
        \frac{\Delta t^{4} }{792}  \bigg(\frac{d^{3}f(t_{n})}{dt^{3}} - \frac{d^{3} {f}(t_{n+1})}{dt^{3}}  \bigg) \nonumber \\
        &+ \frac{\Delta t^{5} }{15840}  \bigg(\frac{d^{4}f(t_{n})}{dt^{4}} + \frac{d^{4} {f}(t_{n+1})}{dt^{4}}  \bigg) +  
        \frac{\Delta t^{6} }{665280}  \bigg(\frac{d^{5}f(t_{n})}{dt^{5}} - \frac{d^{5} {f}(t_{n+1})}{dt^{5}}  \bigg)
        \nonumber \\
        &+ \textbf{J}_{\texttt{H12}}(\Delta  t_{\times}, \Delta t), \ \ \ \ 
    \label{app3_disco_time_h12}
\end{align}
\begin{figure*}
\includegraphics[width=84mm]{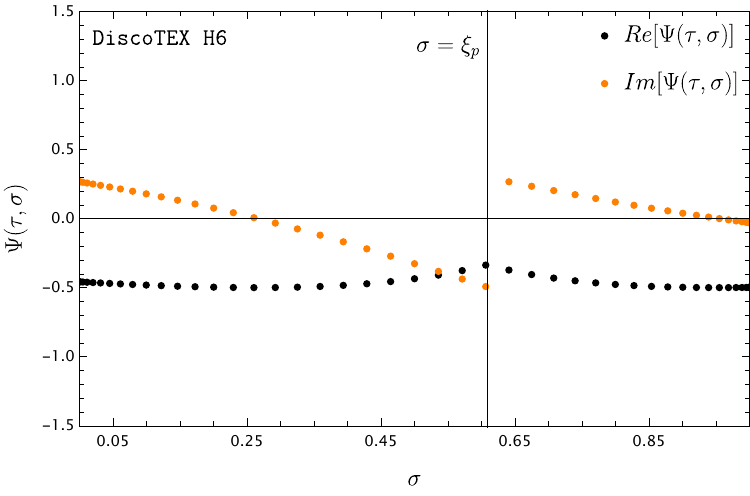}
\quad
\includegraphics[width=84mm]{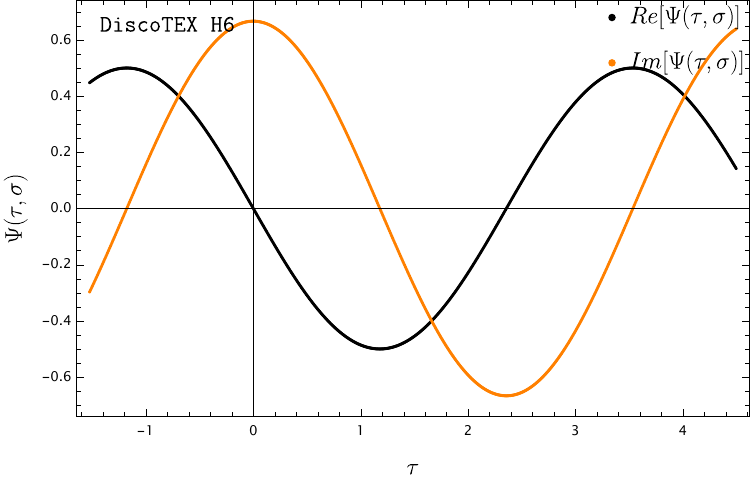}
\caption{Numerical weak-form solution to $\Psi(\tau,\sigma)$ obtained via the \texttt{DiscoTEX H6} \texttt{6th}- order algorithm. \textbf{Left:} Numerical field $\Psi(\tau,\sigma)$ for a point-particle in time-dependent linear motion $\xi_{p}(\tau_{c}) $, where $v$ is the particle's velocity. Specifically, here $\xi_{p} \approx 0.608$ at the coordinate-time $\tau_{c} \approx0.876$ and $v=1/4$. 
\textbf{Right:} Waveform for the point-particle computed on the numerical domain $\sigma \in [0,1]$ and $\tau \in [-1.52, 4.50]$. As expected the result matches that observed for the solution obtained with \texttt{DiscoTEX H4} \texttt{4th}-order numerical algorithm \cite{da2024discotex}. \label{ch3_discotex_hyperboloidal_wave_snap_h6}} 
\end{figure*} 
where $\textbf{J}_{\texttt{H12}}(\Delta  t_{\times}, \Delta t)$ is given by, 
\begin{small}
\begin{align}
  &\textbf{J}_{\texttt{H12}}(\Delta  t_{\times}, \Delta t) = \frac{1}{2}(\Delta t - 2 \Delta  t_{\times})  \textbf{J}_{0}  
   + \frac{1}{44}(5 \Delta t^{2} - 22 \Delta t \Delta  t_{\times} +22 \Delta  t_{\times}^{2})\nonumber \\
    &+\frac{1}{132}(\Delta t - 2 \Delta  t_{\times}) (2\Delta t^{2} - 11 \Delta t \Delta  t_{\times} + 11 \Delta  t_{\times}^{2}) \textbf{J}_{2}  \nonumber \\
    &+ \frac{1}{792}(\Delta t^{4} - 12 \Delta t^{3}\Delta  t_{\times} + 45 \Delta t^{2}\Delta  t_{\times}^{2} - 66 \Delta t \Delta  t_{\times}^{3} + 33 \Delta  t_{\times}^{4})\textbf{J}_{3} \nonumber \\    
    &+ \frac{1}{15840}(\Delta t - 2 \Delta  t_{\times})(\Delta t^{4} - 18 \Delta t^{3} \Delta  t_{\times} \nonumber \\
    &+ 84 \Delta t^{2} \Delta  t_{\times}^{2} - 132 \Delta t^{2} \Delta  t_{\times}^{2} + 66 \Delta  t_{\times}^{4})\textbf{J}_{4} \nonumber \\
    &+ \frac{1}{665280} (\Delta t^{6} - 42 \Delta t^{5} \Delta  t_{\times} + 420 \Delta t^{4} \Delta  t_{\times} ^{2} - 1680 \Delta t^{3} \Delta  t_{\times} ^{3} \nonumber \\
    &+ 3150 \Delta t^{2} \Delta  t_{\times}^{4} - 2772 \Delta t \Delta  t_{\times}^{5} + 924 \Delta  t_{\times}^{5}) \textbf{J}_{5} \nonumber \\
    &+ \frac{1}{30240} (\Delta t - \Delta  t_{\times})\Delta  t_{\times} 
    (\Delta t^{4} - 14 \Delta t^{3} \Delta  t_{\times} + 56 \Delta t^{2} \Delta  t_{\times}^{2} \nonumber \\
    &- 84\Delta t \Delta  t_{\times}^{3} + 42 \Delta  t_{\times}^{4})
    \textbf{J}_{5} \nonumber \\
    &- \frac{1}{665280} (\Delta t - 2\Delta  t_{\times})(\Delta t - \Delta  t_{\times} ) \Delta  t_{\times}\nonumber \\
    & \times  (\Delta t^{4} - 18 \Delta t^{3} \Delta  t_{\times}  + 45 \Delta t^{2} \Delta  t_{\times}^{2} - 132 \Delta t \Delta  t_{\times}^{3} + 66\Delta  t_{\times}^{4}) \textbf{J}_{6} \nonumber \\
    &+ \frac{1}{1330560} (\Delta t - \Delta  t_{\times})^{2}\Delta  t_{\times}^{2} 
    (\Delta t^{4} - 12 \Delta t^{3} \Delta  t_{\times} + 45 \Delta t^{2} \Delta  t_{\times}^{2} \nonumber \\
    &- 66 \Delta t \Delta  t_{\times}^{3} + 33 \Delta  t_{\times}^{4}) \textbf{J}_{7} \nonumber \\
    &+\frac{1}{7983360} (\Delta t - 2\Delta  t_{\times})(\Delta t - \Delta t_{\times})^{3}
    \Delta  t_{\times}^{3} (2 \Delta t^{2}  - 11 \Delta t \Delta  t_{\times}  + 11 \Delta  t_{\times} ^{2}) \textbf{J}_{8}  \nonumber \\
    &+\frac{1}{79833600} (\Delta t - \Delta  t_{\times})^{4}\Delta  t_{\times}^{4} (5 \Delta t^{2} - 22 \Delta t \Delta  t_{\times}  + 22 \Delta  t_{\times}^{2})\textbf{J}_{9} \nonumber \\
   &- \frac{1}{79833600}(\Delta t - 2\Delta  t_{\times}) (\Delta t - \Delta  t_{\times})^{5}\Delta  t_{\times}^{5}  \textbf{J}_{10}  \nonumber \\
   &+ \frac{1}{479001600} (\Delta t - \Delta  t_{\times})^{6}\Delta  t_{\times}^{6}  \textbf{J}_{11},
    \label{app3_disco_time_Jh12}
\end{align}
\end{small}

With the aid of \ref{app_disco_time_higherorder} we study the numerical rate of convergence of the discontinuous Hermite interpolation scheme by computing numerical solutions from \texttt{2nd} to \texttt{12th} order to the following time-dependent Legendre polynomial, 
\begin{equation}
    f(t) = P_{5}(t)\Theta(t) + Q_{5}(t)\Theta(t), 
    \label{ch3_legendreQP_time}
\end{equation}
where $P_{5}(t), Q_{5}(t)$ are the fifth Legendre polynomials of the first and second kind respectively. Introducing a discontinuity at $t_{\times}=0$ in the interval $t \in [-0.55,0.45]$, equation \eqref{ch3_legendreQP_time} admits the following analytical jumps, 
\begin{align}
    &\textbf{J}_{0} = \frac{8}{15},  \ \textbf{J}_{1} = \frac{15}{8}, \ \textbf{J}_{2} = -16, \ \textbf{J}_{3} = -\frac{105}{2}, \\
    &\textbf{J}_{4} = 384, \ \textbf{J}_{5} = 945,  \ \textbf{J}_{6} = -3840, \ \textbf{J}_{7} =0, \\
    &\textbf{J}_{8} = -46080, \ \textbf{J}_{9} = 0,  \ \textbf{J}_{10} = -1935360, \ \textbf{J}_{11} = 0. 
    \label{ch3_legende_anal_jumps}
\end{align}
Integrating analytically the above polynomial in the interval $t \in [-0.55,0.45]$ we get $\approx 0.1125883303464025$. From Figure \ref{ch2_Disco_Time_AllOrders} it is observed that the error associated with the discontinuous time-integration rules 
scales with the respective order of the scheme that is used. To further complement our results we also include a comparison if only a smooth time-integration was used, as described by equations (\eqref{app3_disco_time_h2} - \eqref{app3_disco_time_Jh10}), now corrected in the scheme by incorporating the discontinuous behaviour through the discontinuous time-integration rules given through equations (\eqref{app3_disco_time_Jh2} - \eqref{app3_disco_time_Jh12}). To further investigate the effect increasing the order of the integrator had on numerical results we computed the total simulation running time to compute the numerical data showcased in Figure \ref{ch2_Disco_Time_AllOrders} for all the \texttt{12th} order discontinuous Hermite interpolation schemes that were used. As it is observed in Table \ref{tab_discontinuousTime-allorders} there is a significant increase in the total simulation running time.   \\
\begin{table}
\centering
\begin{small}
\begin{tabular}{||l|| c|| }
\hline \hline 
\textrm{Order of \texttt{NDH}}& 
\textrm{Wall-clock times, s} \\ \hline \hline 
Order 2, \texttt{NDH2} & $4.15227$ \ \  \\ 
Order 4, \texttt{NDH4} & $55.5743$ \ \  \\ 
Order 6, \texttt{NDH6} & $122.499$ \ \  \\ 
Order 8, \texttt{NDH8}  & $199.787$ \ \  \\ 
Order 10, \texttt{NDH10} & $268.464$ \ \  \\ 
Order 12, \texttt{NDH12}  & $348.885$ \ \ \\ \hline \hline 
\end{tabular}
\caption{Total wall-clock computational running time to perform the numerical convergence test summarised on Figure \ref{ch2_Disco_Time_AllOrders}.}\label{tab_discontinuousTime-allorders}
\end{small}
\end{table}

A third user-specifiable numerical optimisation control-factor emerges \cite{da2024discotex, da2023hyperboloidal, phdthesis-lidia}
\begin{enumerate}
    \item \textcolor{gray}{\texttt{[CTRL F\ref{controlfactor1}]} - Number of \texttt{N} nodes;}
    \item \textcolor{gray}{\texttt{[CTRL F\ref{controlfactor2}]} - Number of \texttt{J} jumps;} 
    \item \texttt{[CTRL F3]} - $\Delta \tau$/$\Delta t$, time step-size.\label{controlfactor3}
\end{enumerate}
As we noted in \cite{da2024discotex} the jumps in the time direction, $\textbf{J}$ \textbf{are fixed and determined by the order of the discontinuous time-integration algorithm} as evidenced by equations (\eqref{app3_disco_time_h12}, \eqref{app3_disco_time_Jh12}) and the respective equations for the different order integrators given in \ref{app_disco_time_higherorder}.

\section{Numerical weak-form solutions to distributionally sourced partial differential equations via \texttt{DiscoTEX} - generalisation to higher orders}\label{higher_at_discotex}

\subsection{\texttt{DiscoTEX} at higher-orders: Algorithm description}\label{higher_at_discotex_theory}

\begin{figure}
\includegraphics[width=80mm]{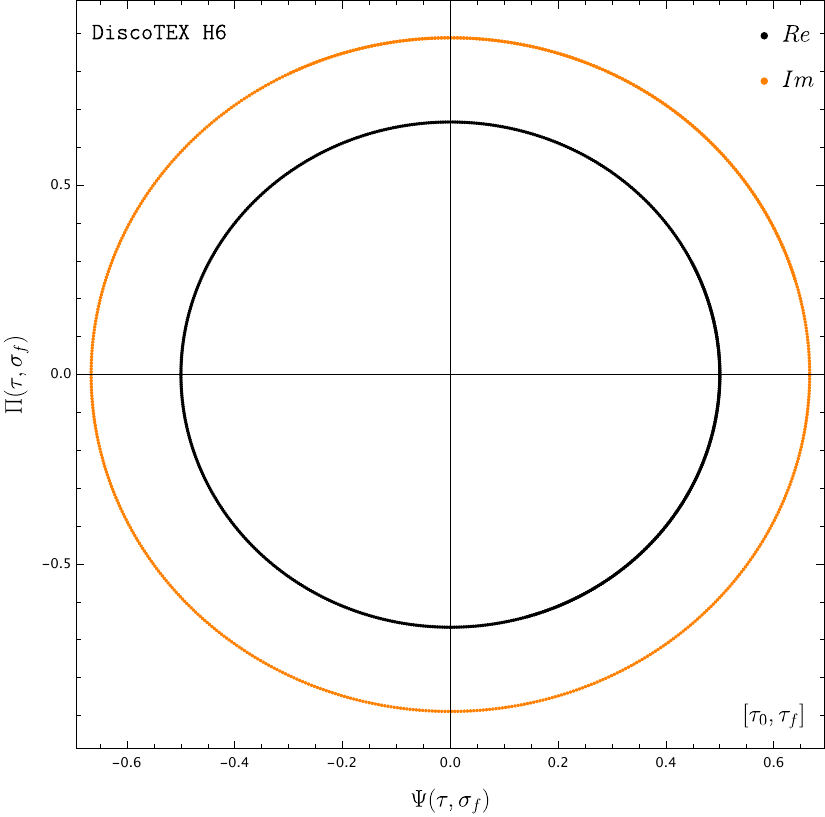}
\caption{Phase portrait for the numerical weak-form solution obtained via \texttt{DiscoTEX H6} with a \texttt{6th}- order \texttt{IMTEX} Hermite \texttt{H6} time integrator in the $(\tau,\sigma)$ hyperboloidal chart. The numerical weak-form solutions to the field $\Psi(\tau,\sigma_{f})$, $\Pi(\tau,\sigma_{f})$ are evaluated in the numerical time interval $\tau \in [-1.52,4.50] $ at the last grid point. As reminded in \cite{da2024discotex} here we too point out that unlike in \cite{da2023hyperboloidal}, Figure 2, uses the exact solutions as initial data, whereas in our previous work, we used trivial initial data, thus not requiring the monitoring of an extra user-specifiable control factor. \label{ch3_discotex_wave_phasepH6}} 
\end{figure}

It will now be demonstrated how \texttt{DiscoTEX} is applied to solve the distributionally-sourced wave-equation of Type I \cite{da2024discotex} given as
\begin{align}
    \square \Psi(\tau,\sigma) = \mathcal{S}(\tau,\sigma) , 
    \label{ch34_hyper_wave_equation}
\end{align}
with, 
\begin{align}
    \square \Psi = \bigg( -\Gamma(\sigma) \partial^{2}_{\tau}  + \varepsilon(\sigma)\partial_{\sigma}\partial_{\tau} + \varrho(\sigma)\partial_{\tau} + \chi(\sigma)\partial^{2}_{\sigma} + \iota(\sigma)\partial_{\sigma}  \bigg) \Psi. 
    \label{hyper_dalembert_op}
\end{align}
The coefficients are specified by
\begin{align}
\label{ch34_sb2_diffOperators_gamma}
&\Gamma(\sigma) = - 4 \sigma^{2}(-1+\sigma^{2}),\\
\label{ch34_sb2_diffOperators_vareps}
&\varepsilon(\sigma) = -4(-1+\sigma) \sigma^{2}(-1 + 2\sigma^{2}),  \\
\label{ch34_sb2_diffOperators_varrho}
&\varrho(\sigma) =  -8(-1+\sigma)\sigma^{3}, \\
\label{ch34_sb2_diffOperators_chi}
&\chi(\sigma) =-4(-1 + \sigma)^{2}\sigma^{4},\\
&\iota(\sigma) = -4 (-1+\sigma)\sigma^{3}(-2+3\sigma),
\label{ch34_sb2_diffOperators_iota}
\end{align}
where, as in \cite{da2024discotex,da2023hyperboloidal}, it's emphasised that the function $\chi(\sigma)$ vanishes at the boundaries where $\sigma = \{0,1\}$, reflecting the outflow behaviour automatically enforced by this hyperboloidal slicing at the future null infinity and the horizon, respectively. The term $\mathcal{S(\tau)}$ is place-hold notation to denote the distributionally source function as only a time-dependent function as before. We will shortly explain how this is handled. As explained in Section \ref{sec_mol} we work via the method-of-lines framework,  
\begin{align}
    &\frac{d\textbf{U}}{d\tau} = \textbf{L} \cdot \textbf{U} + \mathcal{S} + \Upsilon, \ \ \  \textbf{U} = \begin{pmatrix}
        \Psi(\tau) \\ \Pi(\tau)\end{pmatrix},\nonumber \\
    &\mathcal{S} = \begin{pmatrix}  0 \\ \textbf{s}_{\Psi}(\tau) +  \textbf{s}_{\Pi}(\tau)\end{pmatrix},
    \ \ \Upsilon = \begin{pmatrix} \Upsilon_{\Psi}(\tau) \\ \Upsilon_{\Pi}(\tau)  \end{pmatrix}, \label{ch34_red_1ode} \
\end{align}
with the evolution operator $\textbf{L}$ described as, 
\begin{align}
    \label{ch34_hyperboloidal_L1_p}
    &\textbf{L}_{1} = \frac{1}{\Gamma(\sigma)}\bigg( \chi(\sigma)\partial^{2}_{\sigma} + \iota(\sigma)\partial_{\sigma} \bigg), \\
    &\textbf{L}_{2} = \frac{1}{\Gamma(\sigma)}\bigg( \varepsilon(\sigma) \partial_{\sigma} - \varrho(\sigma) \bigg),
    \label{ch34_hyperboloidal_L2_p}
\end{align}
where, for convenience as described in \cite{da2024discotex, da2023hyperboloidal}, the tilde notation is introduced to denote division of the coefficients in the operators by $\Gamma(\sigma)$, e.g., $\tilde{\varepsilon}(\sigma) =\varepsilon(\sigma)/\Gamma(\sigma)$. As motivated in \cite{da2024discotex} we choose to work with the hyperboloidal chart $(\tau,\sigma)$ via the coordinate map $(t \rightarrow t(\tau, \sigma), x \rightarrow x(\sigma))$ described by the \enquote{scri-fixing} technique \cite{zenginouglu2008hyperboloidal, zenginouglu2009gravitational}, 
\begin{align}
    &t = \tau - H(\sigma), \hspace{0.2cm} x = \frac{1}{2} \bigg( \frac{1}{\sigma} + \ln{(1-\sigma)} - \ln{\sigma} \bigg), 
    \label{ch34_minimal_gauge_tx}
\end{align}
where height function $H(\sigma)$ is given by, 
\begin{align}
    &H(\sigma) = \frac{1}{2} \bigg[ \ln{(1-\sigma)} - \frac{1}{\sigma} + \ln{\sigma} \bigg]
    \label{ch34_heightfunction}
\end{align}
as originally introduced by \cite{ansorg2016spectral, jaramillo2021pseudospectrum, macedo2018hyperboloidal, macedo2020hyperboloidal}.

The term $\Upsilon$ is defined as, 
\begin{align}
    \begin{pmatrix} \Upsilon_{\Psi}(\tau) =   [[\Psi]]_{i}\Xi_{i}= J_{0}(\tau)*\Xi_{i}\\ \Upsilon_{\Pi}(\tau) =  [[\Pi]]_{i}\Xi_{i} = \mathbb{J}_{0}(\tau)*\Xi_{i}\end{pmatrix}.
    \label{ch34_td_corrections_tofinal_sol_hyper}
\end{align}
The reader is directed to \ref{app_supplement_disco_spatial}, where the explicit form of the initialising jumps, $J_{0}(\tau), J_{1}(\tau)$ are given. The time jumps are obtained as 
\begin{align}
    \label{ch34_rec_time_PI_jumps_hyper}
    \mathbb{J}(\tau) = \mathbb{J}_{m}(\tau) = \partial_{\tau}(J_{m}(\tau)) - \dot{\xi}_{p}(\tau) J_{m+1}(\tau), \\
    \mathbb{J}(\tau)|_{m=0} = \mathbb{J}_{0}(\tau) = \partial_{\tau}(J_{0}(\tau)) - \dot{\xi}_{p}(\tau) J_{1}(\tau). 
    \label{ch34_rec_time_PI_jump_hyper}
\end{align}
$[[\Psi]]_{i}$ is as defined in equation \eqref{ch2_timedependent_jump}, and $[[\Pi]]_{i}$ is given by 
\begin{align}
    \label{ch3_time-dependent_pi_jump}
    &[[\Pi]]_{i} = \Pi^{+}(\tau,\xi_{p}(\tau)) - \Pi^{-}(\tau,\xi_{p}(\tau)), \hspace{0.2cm} \tau \rightarrow \frac{\sigma_{i}}{v}\\
    &\Xi_{i} = \Theta(\tau_{n+1} - \tau_{i})\Theta(\tau_{i} - \tau_{n}),
    \label{ch34_switch}
\end{align}
where $\Theta$ is as defined in equation \eqref{ch2_HeavisideStepFunction} and $\Xi_{i}$ acts as a \texttt{switch} which turn on these corrections in the time direction for both the $\Psi(\tau)$ and $\Pi(\tau)$ master functions when the particle worldline $x_{p}(t)$ crosses the $i-$th grid-point as they are integrated through \texttt{DiscoTEX} at a time $\tau_{i} \in [\tau_{n}, \tau_{n+1}]: \xi_{p}(\tau)= \sigma_{i}$. We now need to correct the discontinuities in space with the machinery outlined in Section \ref{sec2_disco_space}. Specifically we correct equations (\eqref{s1psi_iota}, \eqref{s2psi_chi}) which contribute to the corrected $\textbf{L}_{1}$ operator associated with the master function $\Psi(\tau,\sigma)$ and equation \eqref{spi_epsilon} contributes to the corrected $\textbf{L}_{2}$ operator. The vector $\mathcal{S}$ in \eqref{ch34_red_1ode} concretely is, 
\begin{align}
    \tilde{\textbf{s}}(\tau) = \begin{pmatrix}
        0 \\
        \tilde{s}_{\Psi,(1)} + \tilde{s}_{\Psi, (2)} + \tilde{s}_{\Pi, (1)}
    \end{pmatrix}, 
    \label{s_spatial_disc_vector_hyper_Wave}
\end{align}

where the subscripts \texttt{(1,2)} denotes \texttt{1st/2nd}- derivative with respect to space. However, as demonstrated in \cite{da2024discotex} higher-order time-derivatives of this source term are necessary, as the orders of the time-integration scheme used are increased. In this work, we use a sixth-order Hermite time-integration rule given as, 
\begin{align}
\textbf{U}_{n+1} - \textbf{U}_{n} &= \int^{\tau_{n+1}}_{\tau_{n}} \bigg( \textbf{L} \cdot \textbf{U} + \mathcal{S} \bigg) d\tau \nonumber \\ 
&= \frac{\Delta \tau}{2} \textbf{L} \cdot (\textbf{U}_{n} + \textbf{U}_{n+1}) + \frac{\Delta \tau}{2} (\mathcal{S}_{n} + \mathcal{S}_{n+1}) \nonumber \\
&+\frac{\Delta \tau ^{2}}{10} \textbf{L} \cdot (\dot{\textbf{U}}_{n} - \dot{\textbf{U}}_{n+1}) + \frac{\Delta \tau^{2}}{10} (\dot{\mathcal{S}}_{n} - \dot{\mathcal{S}}_{n+1}) \nonumber \\ 
&+\frac{\Delta \tau ^{3}}{120} \textbf{L} \cdot (\ddot{\textbf{U}}_{n} + \ddot{\textbf{U}}_{n+1}) + \frac{\Delta \tau^{3}}{120} (\ddot{\mathcal{S}}_{n} + \ddot{\mathcal{S}}_{n+1}) \nonumber \\ 
&= \frac{\Delta \tau}{2} (\textbf{U}_{n} + \textbf{U}_{n+1}) + \frac{\Delta \tau}{2} (\mathcal{S}_{n} + \mathcal{S}_{n+1}) \nonumber \\
&+\frac{\Delta \tau^{2}}{10} \textbf{L} \cdot (\textbf{L} \cdot \textbf{U}_{n} - \textbf{L} \cdot \textbf{U}_{n+1}) + \frac{\Delta \tau^{2}}{10} \textbf{L} \cdot (\mathcal{S}_{n} - \mathcal{S}_{n+1}) \nonumber \\ 
&+\frac{\Delta \tau^{3}}{120} \textbf{L} \cdot (\textbf{L} \cdot \textbf{L} \cdot \textbf{U}_{n} + \textbf{L} \cdot \textbf{L} \cdot \textbf{U}_{n+1}) \nonumber \\
&+ \frac{\Delta \tau^{3}}{120} \textbf{L} \cdot ( \textbf{L} \cdot \mathcal{S}_{n} + \textbf{L} \cdot \mathcal{S}_{n+1})  + \frac{\Delta \tau^{3}}{120} \textbf{L} \cdot (\dot{\mathcal{S}}_{n} + \dot{\mathcal{S}}_{n+1})   \nonumber \\ 
&+ \frac{\Delta \tau^{2}}{10} (\dot{\mathcal{S}}_{n} - \dot{\mathcal{S}}_{n+1}) + \frac{\Delta \tau^{3}}{120} (\ddot{\mathcal{S}}_{n} + \ddot{\mathcal{S}}_{n+1}), 
    \label{discoHermitereduced}
\end{align}
where the following replacement was used, 
\begin{align}
    \textbf{U}^{(2)} &= \frac{d}{d\tau} \bigg[ \textbf{L} \cdot \textbf{U} + \mathcal{S} \bigg], \nonumber \\
    &= \textbf{L} \cdot \big( \textbf{L} \cdot \textbf{U} + \mathcal{S} \big) + \frac{d\mathcal{S}}{d\tau}, \nonumber \\
    &= \textbf{L} \cdot \textbf{L} \cdot \textbf{U} + \textbf{L} \cdot \mathcal{S} + \mathcal{S}^{(1)}, 
    \label{discotex_h6_ddotuvector}
\end{align}
and the superscript \texttt{(2)} denotes the \texttt{2nd}-order derivative with respect to the time coordinate. In \ref{app_higherorder_discotex} we give this necessary transformation for schemes of higher-order i.e. \texttt{8th}- to \texttt{12th}- order, throughout equations (\eqref{discotex_h8_ddotuvector}-\eqref{discotex_h12_ddotuvector}). With a sixth order algorithm the following time-derivative source terms are necessary, 
\begin{align}
    \tilde{\textbf{s}}^{(1)}(\tau) = \begin{pmatrix}
        0 \\
        \tilde{s}^{(1)}_{\Psi,(1)} + \tilde{s}^{(1)}_{\Psi, (2)} + \tilde{s}^{(1)}_{\Pi, (1)}
    \end{pmatrix}, 
    \label{s_spatial_disc_vector_hyper_Wave_dot_h6}
\end{align}
and
\begin{align}
    \tilde{\textbf{s}}^{(2)}(\tau) = \begin{pmatrix}
        0 \\
        \tilde{s}^{(2)}_{\Psi,(2)} + \tilde{s}^{(2)}_{\Psi, (2)} + \tilde{s}^{(2)}_{\Pi, (1)}
    \end{pmatrix}. 
    \label{s_spatial_disc_vector_hyper_Wave_ddot_h6}
\end{align}

\begin{figure*}
    \centering
    \includegraphics[width=16cm]{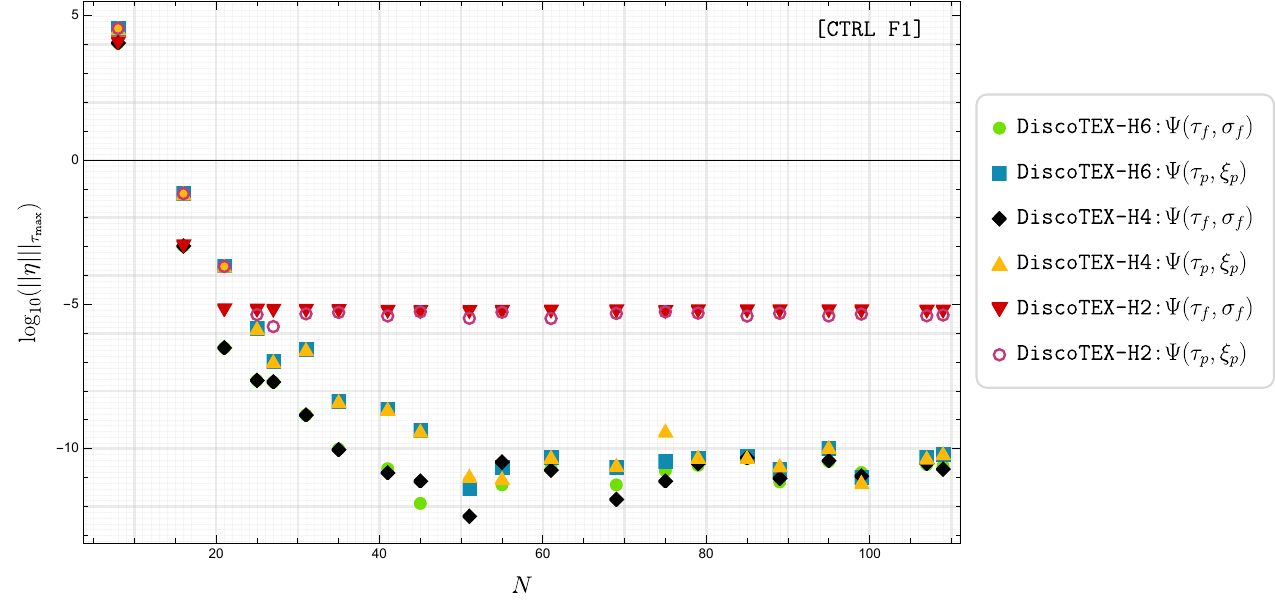}
    \quad 
    \includegraphics[width=16cm]{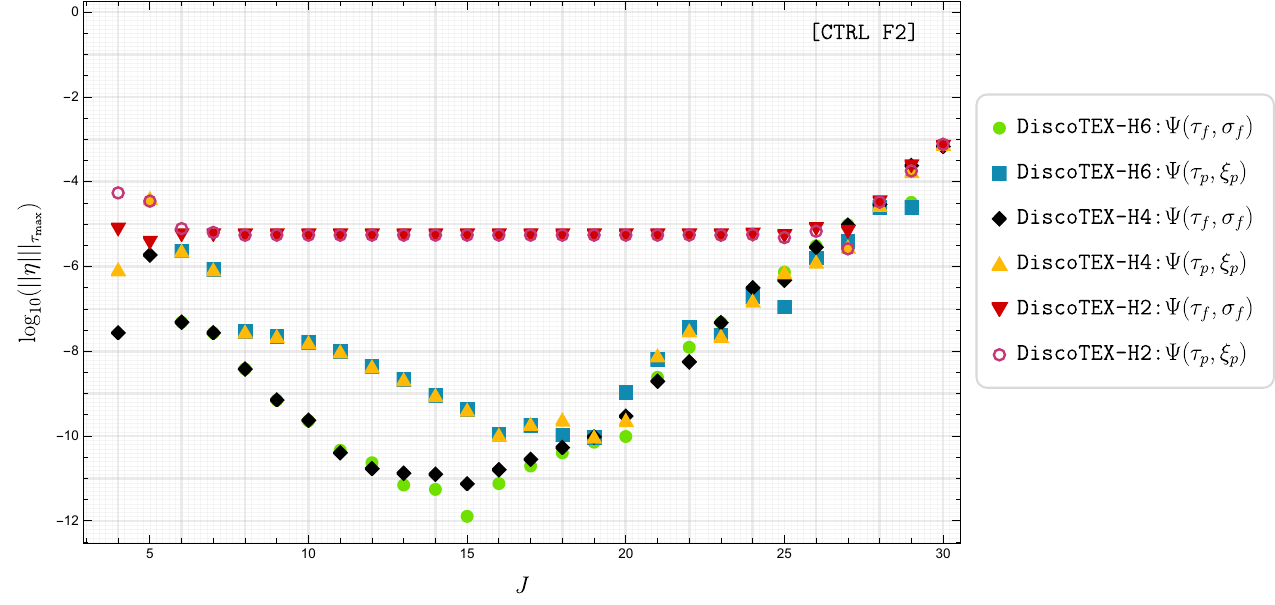}
    \caption{Numerical convergence studies assessing the optimal user-specifiable control factors: number of \texttt{N} nodes - \texttt{[CTRL F\ref{controlfactor1}]} (top plot) and number of \texttt{J} jumps - \texttt{[CTRL F\ref{controlfactor2}]} (bottom plot). Whereas there is a significant accuracy improvement after the second-order time-integration scheme, the difference between the \texttt{4th} and \texttt{6th} time-stepper is marginal for a similar numerical set-up. It is reasonable to pick \texttt{N=45} Chebyshev collocation nodes and \texttt{J=19} jumps. }
    \label{wave_ctrlf1_2}
\end{figure*}

The differential operators corrected with the discontinuous collocation algorithm are given by, 
\begin{align}
  &\tilde{\upchi}(\sigma) \partial_{\sigma}^{2}\Psi\bigg|_{\sigma=\sigma_{i}} \nonumber \\
  &= \sum^{N}_{j=0} \bigg(\tilde{\upchi}(\sigma_{i}) \times D_{(2)} \bigg)_{ij} \big[ \Psi_{j} + \Delta_{\Psi}(\sigma_{j}-\sigma_{p};\sigma_{i} - \sigma_{p}) \nonumber \\ 
  &+ \Delta_{\Psi}^{(1)}(\sigma_{j}-\sigma_{p};\sigma_{i} - \sigma_{p}) + \Delta_{\Psi}^{(2)}(\sigma_{j}-\sigma_{p};\sigma_{i} - \sigma_{p})  \big], \nonumber   \\
  \label{chid2_discretization}
\end{align} 

\begin{align}
  &\tilde{\iota}(\sigma) \partial_{\sigma}\Psi\bigg|_{\sigma=\sigma_{i}} = \nonumber \\
  &\sum^{N}_{j=0} \bigg( \tilde{\iota}(\sigma_{i}) \times D_{(1)} \bigg)_{ij} \big[ \Psi_{j} + \Delta_{\Psi}(\sigma_{j}-\sigma_{p};\sigma_{i} - \sigma_{p}) \nonumber \\ 
  &+ \Delta_{\Psi}^{(1)}(\sigma_{j}-\sigma_{p};\sigma_{i} - \sigma_{p}) + \Delta_{\Psi}^{(2)}(\sigma_{j}-\sigma_{p};\sigma_{i} - \sigma_{p}) \big],  \nonumber   \\
  \label{iotad1_discretization}
\end{align}
\begin{align}
  &\tilde{ \upvarepsilon}(\sigma) \partial_{\sigma}\Pi\bigg|_{\sigma=\sigma_{i}} = \nonumber \\
  &\sum^{N}_{j=0} \bigg( \tilde{\upvarepsilon}(\sigma_{i}) \times D_{(1)} \bigg)_{ij} \big[ \Pi_{j} + \Delta_{\Pi}(\sigma_{j}-\sigma_{p};\sigma_{i} - \sigma_{p}) \nonumber \\ 
  &+ \Delta_{\Pi}^{(1)}(\sigma_{j}-\sigma_{p};\sigma_{i} - \sigma_{p})  + \Delta_{\Pi}^{(2)}(\sigma_{j}-\sigma_{p};\sigma_{i} - \sigma_{p}) \big].  \nonumber   \\
  \label{alpha_d1_discretization}
\end{align}

The explicit form of $\tilde{\textbf{s}}(\tau)$ in equation \eqref{s_spatial_disc_vector_hyper_Wave} containing all the necessary corrections to the differential operators is
\begin{align}
    &\tilde{s}_{\Psi,(1)} = \sum^{N}_{j=0} \bigg( \tilde{\iota}(\sigma_{i}) \times D_{(1)} \bigg)_{ij} \big[\Delta_{\Psi}(\sigma_{j}-\sigma_{p};\sigma_{i} - \sigma_{p}) \big], \ \ \ 
    \label{s1psi_iota}
\end{align}
\begin{align}
    &\tilde{s}_{\Psi, (2)} = \sum^{N}_{j=0} \bigg(\tilde{\upchi}(\sigma_{i}) \times D_{(2)} \bigg)_{ij} \big[ \Delta_{\Psi}(\sigma_{j}-\sigma_{p};\sigma_{i} - \sigma_{p}) \big], \ \ \ 
    \label{s2psi_chi}
\end{align}
\begin{align}
    &\tilde{s}_{\Pi, (1)} = \sum^{N}_{j=0} \bigg( \tilde{ \upvarepsilon}(\sigma_{i}) \times D_{(1)} \bigg)_{ij} \big[  \Delta_{\Pi}(\sigma_{j}-\sigma_{p};\sigma_{i} - \sigma_{p}) \big]. 
    \label{spi_epsilon}
\end{align}
Similarly that of $\tilde{\textbf{s}}^{(1)}(\tau)$ in equation we have, \eqref{s_spatial_disc_vector_hyper_Wave_dot_h6} is 
\begin{align}
    &\tilde{s}_{\Psi,(1)}^{(1)} = \sum^{N}_{j=0} \bigg( \tilde{\iota}(\sigma_{i}) \times D_{(1)} \bigg)_{ij} \big[\Delta_{\Psi}^{(1)}(\sigma_{j}-\sigma_{p};\sigma_{i} - \sigma_{p}) \big], \ \ \ 
    \label{s1psi_iota_dot}
\end{align}
\begin{align}
    &\tilde{s}_{\Psi,(2)}^{(1)} = \sum^{N}_{j=0} \bigg(\tilde{\upchi}(\sigma_{i}) \times D_{(2)} \bigg)_{ij} \big[ \Delta_{\Psi}^{(1)}(\sigma_{j}-\sigma_{p};\sigma_{i} - \sigma_{p}) \big], \ \ \ 
    \label{s2psi_chi_dot}
\end{align}
\begin{align}
    &\tilde{s}_{\Pi,(1)}^{(1)} = \sum^{N}_{j=0} \bigg( \tilde{ \upvarepsilon}(\sigma_{i}) \times D_{(1)} \bigg)_{ij} \big[ \Delta_{\Pi}^{(1)}(\sigma_{j}-\sigma_{p};\sigma_{i} - \sigma_{p}) \big], 
    \label{spi_epsilon_dot}
\end{align}
and finally for $\tilde{\textbf{s}}^{(2)}(\tau)$ in equation \eqref{s_spatial_disc_vector_hyper_Wave_ddot_h6}
\begin{align}
    &\tilde{s}_{\Psi,(1)}^{(2)} = \sum^{N}_{j=0} \bigg( \tilde{\iota}(\sigma_{i}) \times D_{(1)} \bigg)_{ij} \big[\Delta_{\Psi}^{(2)}(\sigma_{j}-\sigma_{p};\sigma_{i} - \sigma_{p}) \big], \ \ \ 
    \label{s1psi_iota_ddot}
\end{align}
\begin{align}
    &\tilde{s}_{\Psi,(2)}^{(2)} = \sum^{N}_{j=0} \bigg(\tilde{\upchi}(\sigma_{i}) \times D_{(2)} \bigg)_{ij} \big[ \Delta_{\Psi}^{(2)}(\sigma_{j}-\sigma_{p};\sigma_{i} - \sigma_{p}) \big], \ \ \ 
    \label{s2psi_chi_ddot}
\end{align}
\begin{align}
    &\tilde{s}_{\Pi,(1)}^{(2)} = \sum^{N}_{j=0} \bigg( \tilde{ \upvarepsilon}(\sigma_{i}) \times D_{(1)} \bigg)_{ij} \big[\Delta_{\Pi}^{(2)}(\sigma_{j}-\sigma_{p};\sigma_{i} - \sigma_{p}) \big].  
    \label{spi_epsilon_ddot}
\end{align}
We further remind the reader as per Section \ref{sec2_disco_space} the $\Delta_{\Psi, \Pi}$ correction are given as, 
\begin{small}
\begin{align}
    &\Delta_{\Psi}(\sigma_{j} - \xi_{p}; \sigma - \sigma_{p}) = \nonumber \\
    &= \bigg[\Theta(\sigma_{i} - \xi_{p}) - \Theta(\sigma_{j} - \xi_{p})\bigg]g_{\Psi}(\sigma_{j} - \xi_{p})  \; \rm{when} \; \sigma = \sigma_{i},  \ \ \ \ \ 
    \label{deltapsi}    
\end{align}
\end{small}
with $g_{\Psi}(\sigma_{j} - \xi_{p})$ defined as 
\begin{equation}
    g_{\Psi}(\sigma_{j} - \xi_{p}) = \sum^{M}_{m=0} \frac{J_{m}}{m!}(\sigma_{j} - \xi_{p})^{m}
    \label{gvec_psi}. 
\end{equation} 
Similarly the evolution operator pertaining to the master field $\Pi(\tau,\sigma)$ is corrected to the discontinuous case via, 
\begin{small}
\begin{align}
    &\Delta_{\Pi}(\sigma_{j} - \xi_{p}; \sigma - \sigma_{p}) = \nonumber \\
    &= \bigg[\Theta(\sigma_{i} - \xi_{p}) - \Theta(\sigma_{j} - \xi_{p})\bigg]g_{\Pi}(\sigma_{j} - \xi_{p})  \; \rm{when} \; \sigma = \sigma_{i}, \ \ \ \ \ 
    \label{deltapi}    
\end{align}
\end{small}
where, 
\begin{equation}
    g_{\Pi}(\sigma_{j} - \xi_{p}) = \sum^{M}_{m=0} \frac{\mathbb{J}_{m}}{m!}(\sigma_{j} - \xi_{p})^{m}
    \label{gevc_pi}, 
\end{equation} 
where the $\mathbb{J}_{m}$ are computed as explained in \eqref{ch34_rec_time_PI_jumps_hyper}. Further as per equations (\eqref{s1psi_iota_dot}-\eqref{spi_epsilon_ddot}) we need the following high-order corrections, 
\begin{small}
\begin{align}
    &\Delta_{\Psi}(\sigma_{j} - \xi_{p}; \sigma - \sigma_{p}) = \nonumber \\
    &= \bigg[\Theta(\sigma_{i} - \xi_{p}) - \Theta(\sigma_{j} - \xi_{p})\bigg]g_{\Psi}^{(1)}(\sigma_{j} - \xi_{p})  \; \rm{when} \; \sigma = \sigma_{i},  \ \ \ \ \ 
    \label{deltapsi}    
\end{align}
\end{small}
where 
\begin{align}
    g_{\Psi}^{(1)}(\sigma- \xi_{p}(\tau)) &= \frac{\partial g_{\Psi}(\sigma-\xi_{p}(\tau))}{\partial \tau} = \nonumber \\ 
    &= \frac{\partial}{\partial \tau} \bigg[ \sum^{M}_{m=0}\frac{J(\tau)}{m!}(\sigma_{j} - \xi_{p} (\tau))^{m}\bigg] \nonumber \\
    &= \sum^{M}_{m=0} \bigg[ \frac{1}{m!}\dot{J}(\sigma_{j} - \xi_{p} )^{m} - \frac{m}{m!} \xi^{(1)}_{p} J (\sigma_{j} - \xi_{p} )^{m-1}  \bigg].  
    \label{gdt}
\end{align}
All that is left to do is to correct $\Delta^{(2)}_{\Psi}$ where the second-order in time $g$ vector is,
\begin{align}
    g^{(2)}(\tau) &= \frac{d^{2}g(\tau)}{d\tau^{2}} = \frac{d}{d\tau}\big(g^{(1)}(\tau) \big) \nonumber \\
   &=\frac{d}{d\tau}\bigg[\frac{1}{m!}J^{(1)}(\tau) w^{m}(\tau) -\frac{m}{m!} x^{(1)}_{p}(\tau) J(\tau) w^{m-1}(\tau)\bigg] \nonumber \\ 
   &= \frac{1}{m!}J^{(2)}(\tau) w^{m}(\tau) - \frac{m}{m!}J(\tau)x^{(2)}_{p}(\tau)w^{m-1}(\tau) \nonumber \\ 
   &+  \frac{m(m-1)}{m!}J(\tau)x^{(1)}_{p}(\tau)^{2}w^{m-2}(\tau) \nonumber \\
   &- 2 x^{(1)}_{p}(\tau) J^{(1)}(\tau) \frac{m}{m!}w^{m-1}(\tau), 
    \label{discotex_h4_ddotg}
\end{align}
\begin{figure*}
   \centering
   \includegraphics[width=16cm]{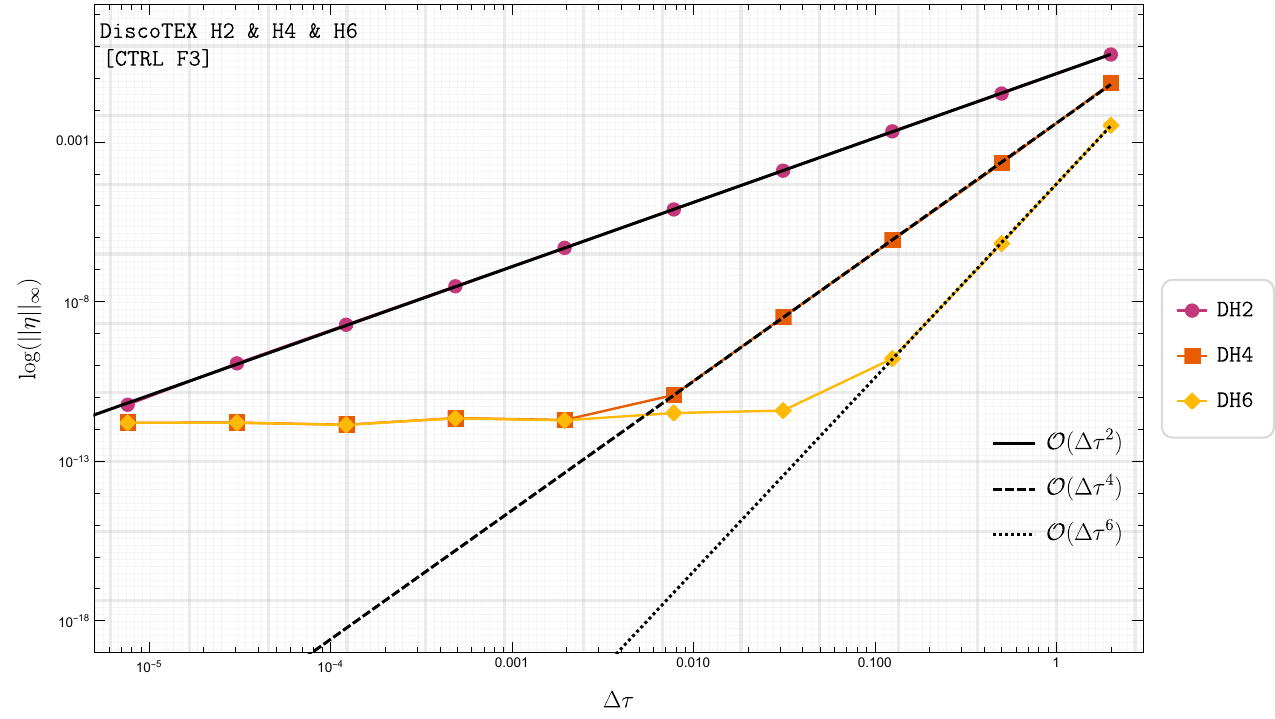}
 \caption{Numerical error associated with computation of the numerical weak-form solution $\Psi(\tau,\sigma_{f})$ against the exact solution in equation \eqref{weakformsolution_wave} with both the discontinuous Hermite integrator of order-2, order-4 and order-6, respectively \texttt{DH2,DH4, DH6}. As expected, we observe in lines \texttt{DH2, DH4, DH6} that the \texttt{DiscoTEX} algorithm converges to its expected orders respectively.}
    \label{wave_ctrlf3}
\end{figure*}
and it trivially follows, for the higher-order time-derivative $\Delta_{\Pi}$ corrections, the jumps used are now $\mathbb{J}_{m}$ as demonstrated in equations (\eqref{deltapi}, \eqref{gevc_pi}). In the \ref{app_higherorder_discotex} we further give the high-order time-derivatives of the $g$ vector needed as the order of the algorithm is further increased in equations (\eqref{discotex_h8_ddotg}, \eqref{discotex_h10_dddotg}, \eqref{discotex_h12_dddotg}). Here, and unlike in \cite{da2024discotex}, we give the corrections without the diagonalisation, and, duly, the dot product has also been removed. As stressed in that work, as long as the operations are handled with care there will be no need to diagonalise these quantities. As longs as we ensure: $1)$ the derivatives are always acting on the governing evolution field - from the left side; and, $2)$ if no differential operators are acting on the field variables on $\textbf{L}_{1/2}$ - like for example the last terms in both equations \eqref{ch34_hyperboloidal_L1_p} and \eqref{ch34_hyperboloidal_L2_p} - we multiply these terms by the identity matrix, for example, $\varrho(\sigma)$ in equation \eqref{ch34_hyperboloidal_L1_p} would be coded as $\varrho(\sigma_{i}) \times \textbf{I}$.  \\

\texttt{DiscoTEX H6}'s numerical algorithm is thus succinctly
\begin{align}
    &\textbf{U}_{n+1} = \textbf{U}_{n}+  \textbf{HFH6} \cdot  \bigg[ \textbf{A} \cdot  \bigg[ \textbf{TEXH6} \cdot \textbf{U}_{n} + \frac{ \Delta \tau}{10} \bigg(\textbf{s}_{n} - \textbf{s}_{n+1} \bigg) \nonumber \\
    &+\frac{ \Delta \tau^{2}}{120} \bigg( \big(\textbf{s}^{(1)}_{n} - \textbf{s}^{(1)}_{n+1} \big) + \textbf{L} \cdot \big(\textbf{s}^{(1)}_{n} - \textbf{s}^{(1)}_{n+1} \big)\bigg) \bigg]\nonumber \\
    &+  \frac{\Delta \tau}{2} \bigg( \textbf{s}_{n} + \textbf{s}_{n+1} \bigg) + \frac{\Delta \tau^{2}}{10} \bigg(\textbf{s}^{(1)}_{n} - \textbf{s}^{(1)}_{n+1}  \bigg) + \frac{\Delta \tau^{3}}{120} \bigg(\textbf{s}^{(2)}_{n} + \textbf{s}^{(2)}_{n+1}  \bigg) \nonumber \\
    &+ \Upsilon +  \textbf{J}_{\texttt{H6}} (\Delta \tau_{\times}, \Delta \tau) \; \Xi \bigg],  
    \label{ch3_discoTEX6_Wave_hyper}
\end{align} 
where, 
\begin{align}
    \textbf{TEXH6} = \bigg[ \textbf{I} + \frac{\textbf{A}\cdot\textbf{A}}{60}\bigg]. 
    \label{discotex12_tex6}
\end{align}
The higher-order \texttt{DiscoTEX H8-H12} numerical algorithms can be found in \ref{app_higherorder_discotex} through equations (\eqref{ch3_discoTEX2_Wave_hyper}-\eqref{ch3_discoTEX12_Wave_hyper}). It is left to specify the jumps associated with the implementation of the higher-order time-integration scheme as described in Section \ref{sec2_disco_time} highlighted in the algorithm by the term $\textbf{J}_{\texttt{H6}} (\Delta \tau_{\times}, \Delta \tau)$. We need 6 jumps for each of the differential operations in the evolution operators (\eqref{ch34_hyperboloidal_L1_p}, \eqref{ch34_hyperboloidal_L2_p}), i.e.
\begin{align}
    \textbf{J}_j(\tau) &= 
    \begin{pmatrix}
    \mathbb{JJ}_j(\tau)\\
    \mathcal{JJ}_j(\tau)
\end{pmatrix} =  \nonumber \\
    &=\begin{pmatrix}
    {\mathbb{J}_{0}(\tau),\mathbb{K}_{0}(\tau),\mathbb{L}_{0}(\tau),\mathbb{M}_{0}(\tau), \mathbb{N}_{0}(\tau), \mathbb{O}_{0}(\tau)} \\
    {\mathcal{J}_{0}(\tau), \mathcal{K}_{0}(\tau),\mathcal{L}_{0}(t),\mathcal{M}_{0}(\tau), \mathcal{N}_{0}(\tau), \mathcal{O}_{0}(\tau)}
\end{pmatrix}\bigg|_{\tau \rightarrow \frac{\sigma_{p}(\tau)}{v}}. 
    \label{discotexh6_timejumps_hyper}
\end{align}
We refer the reader to \ref{app_hyper_timejumps} where the explicit expressions for $\textbf{J}_{0}(\tau)$ to $\textbf{J}_{5}(\tau)$ are given through equations (\eqref{hyper_wave_jb_0}-\eqref{wave_Mc_0_exp}). Additionally for higher-order \texttt{DiscoTEX}'s implementations the following jumps are necessary: - \texttt{DiscoTEX H8}'s will require up to $\textbf{J}_{7}(\tau)$ jumps; - \texttt{DiscoTEX H10}'s will require up to $\textbf{J}_{9}(\tau)$ jumps, and finally \texttt{DiscoTEX H12}'s will require up to $\textbf{J}_{11}(\tau)$ jumps. 

\subsection{\texttt{DiscoTEX} at higher-orders: Numerics discussion}\label{higher_at_discotex_numerics}
\begin{table*}
\centering
\begin{small}
\begin{tabular}{||l|| c  || c ||c || c|| }
\hline \hline 
\texttt{Algorithm}&
$\eta$ &
\texttt{Wall-clock times, s}&
\texttt{[\# CTRL Factors]}\\ \hline \hline 
\rowcolor{retroorange}\texttt{DiscoTEX H2}& $5.5\times10^{-6}$ &10.8018268& 3 \\ \hline \hline 
\rowcolor{mygreen}\texttt{DiscoTEX H4}&$7.7\times10^{-11}$&11.1215158 & 3\\ \hline
\rowcolor{mygreen}\texttt{DiscoTEX H6} & $7.2\times10^{-11}$ &42.847469& 3 \\ \hline \hline 
\rowcolor{retroorange}\texttt{DiscoTEX H8}&$7.2\times10^{-9}$& 147.7311478 & 3\\ 
\rowcolor{retroorange}\texttt{DiscoTEX H10} & $7.2\times10^{-11}$ & 517.8653133& 3\\
\rowcolor{retroorange}\texttt{DiscoTEX H12}&$7.2\times10^{-11}$& 1570.6355016& 3\\ \hline \hline 
\end{tabular}
\end{small}
\caption{Comparison between \texttt{DiscoTEX H2-H12}. As hinted by Table 8 of reference \cite{da2024discotex} increasing the order of integration has not resulted in a significant accuracy improvement for a given time step that would justify the extra computational cost. \texttt{DiscoTEX H4} offers the best compromise with accuracy and minimal simulation running time with the current numerical framework.} 
\label{table_higher_at_discotex} 
\end{table*}
In this work, the focus is on demonstrating the machinery necessary for implementing the \texttt{DiscoTEX} algorithm at higher orders of numerical convergence and how the performance compares. In Figures \ref{ch3_discotex_hyperboloidal_wave_snap_h6} and \ref{ch3_discotex_wave_phasepH6} the behaviour is the same as for \texttt{DiscoTEX H4}'s Figures 8 and 9 \cite{da2024discotex}, as expected. We ran our simulations in the same set-up as in the companion paper described in Section 3.4.3, that is the algorithm spans a temporal domain of $\tau \in [-1.52, 4.50]$ and a spatial domain of $\sigma = [0,1]$, where, as expected from the minimal gauge implementation future null infinity $\mathscr{I}^{+}$ is mapped to $\sigma=0$ and the horizon to $\sigma=1$. The relative error is computed as
\begin{equation}
    \eta = \bigg| 1- \frac{\Psi_{\rm{Numerical}}}{\Psi_{\rm{Exact}}} \bigg|.
    \label{ch3_relativeerror}
\end{equation}

As reviewed in Section \ref{discotex_nutshell} three numerical optimisation control factors need to be studied to ensure optimal implementation of the \texttt{DiscoTEX} algorithm: \texttt{[CTRL F\ref{controlfactor1}]} - number of \texttt{N} nodes; \texttt{[CTRL F\ref{controlfactor2}]} - number of \texttt{J} jumps and   \texttt{[CTRL F\ref{controlfactor3}]} - timestep size. Following up from Section \ref{sec2_disco_time}, and to verify whether the algorithm has been implemented correctly, we determine \texttt{[CTRL F\ref{controlfactor3}]} and check numerical convergence order is as expected for the given algorithms under study. In Figure \ref{wave_ctrlf3} there are no surprises and we observe \texttt{DiscoTEX H2, H4 \& H6} converge to their respective orders.

We then assessed \texttt{[CTRL F\ref{controlfactor1}]} and \texttt{[CTRL F\ref{controlfactor2}]} and determined \texttt{N=45} and \texttt{J=19} to result in an optimal implementation of \texttt{DiscoTEX H6}. Informed by Figure \ref{controlfactor3}, and to ensure fairness in the comparison provided between the different order algorithmic implementations, we select for the numerical simulations the same timestep used in \cite{da2024discotex} i.e. $\Delta \tau = 0.00666667$. For the assessment of \texttt{[CTRL F\ref{controlfactor1}]}, 15 jumps were used whereas the numerical tests for determining \texttt{[CTRL F\ref{controlfactor2}]} were computed with \texttt{N+1=45} Chebyshev-Gauss-Lobatto collocation nodes. From Figure \ref{wave_ctrlf1_2} we can see that the results for the higher-order algorithm closely follow the numerical convergence pattern observed for \texttt{DiscoTEX H4} without any substantial improvements. To further complement our results we investigate the computational times for all 12 orders of the \texttt{DiscoTEX} algorithm. From Table \ref{table_higher_at_discotex} it is very clear that, for the same numerical control factors, after the sixth-order algorithm the accuracy stagnates while the computational time significantly increases with increasing order of the time-stepper. Compared with Table 8 of the companion paper \cite{da2024discotex}, which computed the simulations for a much longer time i.e. up to $\tau =10,000$ amounting to half a million timesteps, the computational times are approximately as long as those observed here which only computes about 903 timesteps. The reason for this significant difference is because of the extra time jumps added by the fact our trajectory is time-dependent and hence one needs to account for the jumps in the time direction, unlike the case studied in that section which probed circular orbits i.e. at a fixed particle position of $\sigma= 0.25$, for example. \\

We finally note in this paper, that, unlike in the companion paper \cite{da2024discotex}, we don't explicitly show any checks to assess whether the algorithm is capable of giving the solution at either limit of the point-particle position. The reason for this \enquote{omission} is that it is a trivial calculation and the procedure is exactly as described in the first paper. In assessing the numerical optimisation \texttt{[CTRL F\ref{controlfactor1} \& F\ref{controlfactor2}]} we also included the result for the interpolation at the particle value, this is demonstrated in Figure \ref{wave_ctrlf1_2}.

\section{Conclusion}\label{finale}

This paper extends the work presented in ref. \cite{da2024discotex} to higher orders of numerical convergence. The main conclusion of this work is the significant increase in the computational time and accuracy stagnation as we increase past \texttt{6th}-order. In this work, we give all the necessary machinery required to implement \newline \texttt{DiscoTEX} including all the time jumps that are to be expected from any point-particle trajectory that has a time dependence. We observed that this caused a significant increase in the total numerical simulation running time for a (much) shorter time interval. \\

For cases where time-dependent point-particle trajectories are unavoidable, such as when modelling highly eccentric EMRIs, higher-order versions of the algorithm should be avoided, especially if conservative self-force computations are to be computed \cite{da2024discotex, da2023hyperboloidal, afshordi2023waveform}. However, naturally as expected, and it is observed in Figures \ref{ch2_Disco_Time_AllOrders} and \ref{wave_ctrlf3}, higher-order implementations tend to reach higher accuracies faster for longer timesteps. This faster convergence could be very useful for other astrophysical modelling scenarios including the modelling of extremely asymmetric binaries such as XMRIs. On top of not requiring any higher-order corrections to the self-force, some XMRIs, depending on the compact object companion, may also more likely be prescribing a circular orbit as they (very) slowly inspiral towards the supermassive black hole \cite{amaro2019extremely, gourgoulhon2019gravitational,  barsanti2022detecting, vazquez2022revised, vazquez2023sgr}. Another area where it may be very useful to use the discontinuous or (smooth!) \texttt{IMTEX} Hermite time-stepper could be black-hole spectroscopy. These applications are left for future work and are to be compared with different numerical algorithms. 

\section*{Acknowledgements}
I would like to thank my PhD supervisors Professors Pau Figueras and Juan Valiente Kroon for their encouragement and guidance in the completion of this work. Finally, it is a pleasure to thank Subrahmanyan Chandrasekhar, whose intelligence and work ethic have been a source of motivation and (inspiration) in the last years of my PhD \cite{chandra, phdthesis-lidia}. 

\appendix

\section{Complementary work for the development of the \texttt{DiscoTEX} numerical algorithm at higher-orders of numerical convergence}

\subsection{Higher-order discontinuous time-integration \newline schemes }\label{app_disco_time_higherorder}

As a necessary complement to Section \ref{sec2_disco_time} the following discontinuous time-integration rules from orders \texttt{2nd} to \texttt{12th} are explicitly given below: 
\begin{align}
  &f(t)_{\texttt{DH2}} = \frac{\Delta t }{2} \bigg(f(t_{n}) + f(t_{n+1} )  \bigg) + \textbf{J}_{\texttt{H2}}(\Delta  t_{\times}, \Delta t), 
\label{app3_disco_time_h2}
\end{align}

\begin{align}
    &f(t)_{\texttt{DH4}} =\frac{\Delta t }{2} \bigg(f(t_{n}) + f(t_{n+1} )  \bigg) + \frac{\Delta t^{2} }{12}  \bigg(\dot{f}(t_{n}) - \dot{f}(t_{n+1} )  \bigg) \nonumber \\
    &+ \textbf{J}_{\texttt{H4}}(\Delta  t_{\times}, \Delta t), 
    \label{app3_disco_time_h4}
\end{align}

\begin{align}
    &f(t)_{\texttt{DH8}} = \frac{\Delta t }{2} \bigg(f(t_{n}) + f(t_{n+1} )  \bigg) + \frac{3}{28} \Delta t^{2} \bigg(\dot{f}(t_{n}) - \dot{f}(t_{n+1} )  \bigg) \nonumber \\
    &+ \frac{\Delta t^{3} }{84}  \bigg(\ddot{f}(t_{n}) + \ddot{f}(t_{n+1} )  \bigg) + 
    \frac{\Delta t^{4} }{1680}  \bigg(\frac{d^{3}f(t_{n})}{dt^{3}} - \frac{d^{3} {f}(t_{n+1})}{dt^{3}}  \bigg) \nonumber \\
    &+ \textbf{J}_{\texttt{H8}}(\Delta  t_{\times}, \Delta t), \ \ \ \ 
    \label{app3_disco_time_h8}
\end{align}

\begin{small}
\begin{align}
&f(t)_{\texttt{NDH6}} = \frac{\Delta t }{2} \bigg(f(t_{\nu}) + f(t_{\nu+1} )  \bigg) + \nonumber \\
&\frac{\Delta t^{2} }{10}  \bigg(\dot{f}(t_{\nu}) - \dot{f}(t_{\nu+1} )  \bigg) + \frac{\Delta t^{3} }{120}  \bigg(\ddot{f}(t_{\nu}) + \ddot{f}(t_{\nu+1} )  \bigg) + \textbf{J}_{\texttt{H6}}(\Delta t_{\times}, \Delta t), \ \ \ \ \ \ \ \ \label{sec2_disco_time_h6}
\end{align}
\end{small}
where $\textbf{J}_{\texttt{H6}}(\Delta t_{\times}, \Delta t)$ is given by
\begin{align}
    &\textbf{J}_{\texttt{H6}}(\Delta t_{\times} , \Delta t) = \frac{1}{2}(\Delta t - 2 \Delta t_{\times} ) \textbf{J}_{0}  + \frac{1}{10} (\Delta t^{2} - 5 \Delta t \Delta t_{\times} + 5 \Delta t_{\times}^{2}) \textbf{J}_{1} + \nonumber\\
    &\frac{1}{120} (\Delta t - 2 \Delta t_{\times}) (\Delta t^{2} - 10 \Delta t \Delta t_{\times}  + 10 \Delta t_{\times} ^{2}) \textbf{J}_{2}  -  \nonumber \\
    &\frac{1}{120}(\Delta t - \Delta t_{\times}) \Delta t_{\times} (\Delta t^{2} - 5\Delta t \Delta t_{\times} + 5 \Delta t_{\times} ^{2}) \textbf{J}_{3} + \nonumber \\
    &\frac{1}{240} (\Delta t - 2\Delta t_{\times} ) (\Delta t - \Delta t_{\times} )^{2} \Delta t_{\times} ^{2}  \textbf{J}_{4} - \nonumber \\
    &\frac{1}{720} (\Delta t - \Delta t_{\times} )^{3}\Delta t_{\times}^{3}   \textbf{J}_{5}. \label{sec2_disco_time_Jh6}
\end{align}

\begin{align}
        &f(t)_{\texttt{DH10}} =\frac{\Delta t }{2} \bigg(f(t_{n}) + f(t_{n+1} )  \bigg) + \frac{\Delta t^{2} }{9}  \bigg(\dot{f}(t_{n}) - \dot{f}(t_{n+1} )  \bigg) \nonumber \\
        &+ \frac{\Delta t^{3} }{72}  \bigg(\ddot{f}(t_{n}) + \ddot{f}(t_{n+1} )  \bigg) + 
        \frac{\Delta t^{4} }{1008}  \bigg(\frac{d^{3}f(t_{n})}{dt^{3}} - \frac{d^{3} {f}(t_{n+1})}{dt^{3}}  \bigg) \nonumber \\
        &+ \frac{\Delta t^{5} }{30240}  \bigg(\frac{d^{4}f(t_{n})}{dt^{4}} + \frac{d^{4} {f}(t_{n+1})}{dt^{4}}  \bigg) \nonumber \\
        &+ \textbf{J}_{\texttt{H10}}(\Delta  t_{\times}, \Delta t), \ \ \ \ 
    \label{app3_disco_time_h10}
\end{align}

where $\textbf{J}_{\texttt{H2-H12}}(\Delta  t_{\times}, \Delta t)$ are given respectively as,  
\begin{align}
    &\textbf{J}_{\texttt{H2}}(\Delta  t_{\times}, \Delta t) 
    = \frac{1}{2}(\Delta t - 2 \Delta  t_{\times}) \textbf{J}_{0} 
    + \frac{\Delta  t_{\times}}{2}(\Delta  t_{\times} - \Delta t) \textbf{J}_{1}, 
    \label{app3_disco_time_Jh2}
\end{align}

\begin{align}
    &\textbf{J}_{\texttt{H4}}(\Delta  t_{\times}, \Delta t) = \frac{1}{2}(\Delta t - 2 \Delta t_{\times} )  \textbf{J}_{0}  + \frac{\Delta  t_{\times}^{2}}{12}(\Delta t^{2} - 6 \Delta t \Delta  t_{\times} + 6 \Delta t_{\times} ^{2})  \textbf{J}_{1} \nonumber \\
    &- \frac{1}{12}\Delta  t_{\times} (\Delta t^{2} - 3 \Delta t \Delta  t_{\times} + 2 \Delta  t_{\times} ^{2}) \textbf{J}_{2} + \frac{1}{24}\Delta  t_{\times}^{2}, (\Delta t-\Delta  t_{\times}) \textbf{J}_{3} \ \ \ \ \ \ \ \ 
    \label{app3_disco_time_Jh4}
\end{align}

\begin{align}
    &\textbf{J}_{\texttt{H8}}(\Delta  t_{\times}, \Delta t) =\frac{1}{2}(\Delta t - 2 \Delta  t_{\times})  \textbf{J}_{0}  + \frac{1}{28} (3\Delta t^{2} - 14 \Delta t \Delta  t_{\times} + 14 \Delta  t_{\times}^{2})  \textbf{J}_{1} + \nonumber \\
    &+\frac{1}{84}(\Delta t - 2 \Delta  t_{\times}) (\Delta t^{2} - 7 \Delta t \Delta t_{\times} + 7 \Delta  t_{\times}^{2}) \textbf{J}_{2}  \nonumber \\
    &+ \frac{1}{1680}(\Delta t^{4} - 20 \Delta t^{3}\Delta  t_{\times} + 90 \Delta t^{2}\Delta  t_{\times}^{2} - 140 \Delta t \Delta  t_{\times}^{3} + 70 \Delta  t_{\times} ^{4})\textbf{J}_{3} \nonumber \\
    &- \frac{1}{1680}(\Delta t - 2 \Delta  t_{\times} )(\Delta t - \Delta  t_{\times} )\Delta  t_{\times}  (\Delta t^{2} - 7\Delta t \Delta  t_{\times} + 7 \Delta  t_{\times}^{2})\textbf{J}_{4} \nonumber \\
    &+ \frac{1}{10080} (\Delta t - \Delta  t_{\times})^{2} \Delta  t_{\times} ^{2}(3 \Delta t^{2} - 14\Delta t \Delta  t_{\times}  + 14 \Delta  t_{\times} ^{2}) \textbf{J}_{5}
    \nonumber \\
    &- \frac{1}{10080} (\Delta t - 2\Delta  t_{\times})(\Delta t - \Delta  t_{\times} )^{3} \Delta  t_{\times}^{3} \textbf{J}_{6} 
   + \frac{1}{40320} (\Delta t - \Delta  t_{\times} )^{4}\Delta  t_{\times} ^{4}  \textbf{J}_{7}, \hspace{12mm} \label{app3_disco_time_Jh8}
\end{align}

\begin{small}
\begin{align}
   &\textbf{J}_{\texttt{H10}}(\Delta  t_{\times}, \Delta t) = \frac{1}{2}(\Delta t - 2 \Delta  t_{\times})  \textbf{J}_{0}  
   + \frac{1}{18}(\Delta t - 3 \Delta  t_{\times}) (2 \Delta t - 3 \Delta  t_{\times}) \nonumber \\
    &+\frac{1}{72}(\Delta t - 2 \Delta  t_{\times}) (\Delta t^{2} - 6 \Delta t \Delta  t_{\times} + 6 \Delta  t_{\times}^{2}) \textbf{J}_{2}  \nonumber \\
    &+ \frac{1}{1008}(\Delta t^{4} - 14 \Delta t^{3}\Delta  t_{\times} + 56 \Delta t^{2}\Delta  t_{\times}^{2} - 84 \Delta t \Delta  t_{\times}^{3} + 42 \Delta  t_{\times}^{4})\textbf{J}_{3} \nonumber \\
    &- \frac{1}{30240}(\Delta t - 2 \Delta  t_{\times})(\Delta t^{4} - 28 \Delta t^{3} \Delta  t_{\times}  + 154 \Delta t^{2} \Delta  t_{\times}^{2} - 252 
    \Delta t \Delta  t_{\times} ^{3} + 126 \Delta  t_{\times}^{4})\textbf{J}_{4} \nonumber \\
    &+ \frac{1}{30240} (\Delta t - \Delta  t_{\times} )\Delta  t_{\times}  
    (\Delta t^{4} - 14 \Delta t^{3} \Delta  t_{\times}  + 56 \Delta t^{2} \Delta  t_{\times}^{2} - 84\Delta t \Delta  t_{\times} ^{3} + 42 \Delta  t_{\times}^{4})
    \textbf{J}_{5} \nonumber \\
    &+ \frac{1}{60480} (\Delta t - 2\Delta  t_{\times})(\Delta t - \Delta  t_{\times})^{2} \Delta  t_{\times} ^{2} (\Delta t^{2} - 6 \Delta t \Delta  t_{\times} + 6 \Delta  t_{\times}^{2}) \textbf{J}_{6} \nonumber \\
    &- \frac{1}{362880} (\Delta t - 3\Delta  t_{\times})(2\Delta t - 3\Delta  t_{\times})(\Delta t - \Delta  t_{\times} )^{3} \Delta  t_{\times}^{3} \textbf{J}_{7} \nonumber \\
    &+ \frac{1}{725760} (\Delta t - 2\Delta  t_{\times})(\Delta t - \Delta  t_{\times})^{4} \Delta  t_{\times}^{4} \textbf{J}_{8} 
   - \frac{1}{3628800} (\Delta t - \Delta  t_{\times})^{5}\Delta  t_{\times}^{5}  \textbf{J}_{9}, \hspace{2mm} 
    \label{app3_disco_time_Jh10}
\end{align}
\end{small}

\subsection{Supplementary material to the spatial discontinuous discretization}\label{app_supplement_disco_spatial}

As showed in \cite{da2024discotex} for the equation in \eqref{ch34_red_1ode} we have the following initialising jumps, 
\begin{align}
    \label{hyper_j0}
    &J_{0} = \bar{\gamma}^{2} \frac{\triangle_{c}(\xi_{p})}{|\triangle_{c}(\xi_{p})|^{3}} F(\tau_{c})   \\
    &J_{1} = \bar{\gamma}^{2} \bigg[  \frac{\triangle_{c}'(\xi_{p})}{|\triangle_{c}(\xi_{p})|^{3}}F(\tau_{c})  - \frac{\triangle_{c}(\xi_{p})}{|\triangle_{c}(\xi_{p})|^{3}} \frac{\partial F(\tau_{c})}{\partial t} \frac{\partial t}{\partial \xi_{p}} +  \frac{G(\tau_{c})}{\triangle_{c}(\xi_{p})}\nonumber \\
    &+ \bigg( \ddot{\xi}_{p}\Gamma(\xi_{p})  + \dot{\xi}_{p}^{2}\Gamma'(\xi_{p}) - \dot{\xi}_{p}\varepsilon'(\xi_{p})  + \varrho'(\xi_{p})   \chi'(\xi_{p}) - \iota(\xi_{p}) \bigg)J_{0} \nonumber \\
    &+ \bigg( 2 \dot{\xi}_{p}\Gamma(\xi_{p})  - \varepsilon(\xi_{p}) \bigg)\dot{J}_{0}  \bigg]. 
    \label{hyper_j1}
\end{align}
where $\tau_{c} =t(\tau,\sigma_{p}(\tau))$, $\xi_{p} = x(\sigma_{p}(\tau))$.  

\subsection{Higher-order \texttt{DiscoTEX} algorithm}\label{app_higherorder_discotex}
For completion, we include all the \texttt{DiscoTEX} algorithms from \texttt{2th}- to \texttt{12th}- order of time-integration. 

The simplest algorithm is \texttt{DiscoTEX H2} given as, 
\begin{align}
    &\textbf{U}_{n+1} = \textbf{U}_{n}+  \textbf{HFH2} \cdot  \bigg[ \textbf{A}  \cdot \textbf{U}_{n} +\frac{\Delta \tau}{2} \bigg(\textbf{s}_{n} + \textbf{s}_{n+1}\bigg) \nonumber \\
    & + \Upsilon +  \textbf{J}_{\texttt{H2}} (\Delta \tau_{\times}, \Delta \tau) \; \Xi \bigg]. 
    \label{ch3_discoTEX2_Wave_hyper}
\end{align} 

\texttt{DiscoTEX H8} is given by,  
\begin{align}
    &\textbf{U}_{n+1} = \textbf{U}_{n}+  \textbf{HFH8} \cdot  \bigg[ \textbf{A} \cdot  \bigg[ \textbf{TEXH8} \cdot \textbf{U}_{n} + \frac{ 3\Delta \tau}{28} \bigg(\textbf{s}_{n} - \textbf{s}_{n+1} \bigg) \nonumber \\
        &+\frac{\Delta \tau^{2}}{84} \bigg( \big(\textbf{s}^{(1)}_{n} +\textbf{s}^{(1)}_{n+1}\big) + \textbf{L} \cdot \big(\textbf{s}_{n} +\textbf{s}_{n+1}\big) \bigg) + \nonumber\\
        &+\frac{\Delta \tau^{3}}{1680} \bigg(\big(\textbf{s}^{(2)}_{n} -\textbf{s}^{(2)}_{n+1}\big) + \textbf{L} \cdot \big(\textbf{s}^{(1)}_{n} -\textbf{s}^{(1)}_{n+1}\big) + \textbf{L} \cdot \textbf{L} \cdot \big(\textbf{s}_{n} -\textbf{s}_{n+1}\big) \bigg) \bigg]\nonumber \\
        &+  \frac{\Delta \tau}{2} \bigg( \textbf{s}_{n} + \textbf{s}_{n+1} \bigg) + \frac{3\Delta \tau^{2}}{28} \bigg(\textbf{s}^{(1)}_{n} -\textbf{s}^{(1)}_{n+1}  \bigg) + \frac{\Delta \tau^{3}}{84} \bigg(\textbf{s}^{(2)}_{n} + \textbf{s}^{(2)}_{n+1}  \bigg) \nonumber \\
        &+ \frac{\Delta \tau^{3}}{1680} \bigg(\textbf{s}^{(3)}_{n} - \textbf{s}^{(3)}_{n+1}  \bigg) + \Upsilon +  \textbf{J}_{\texttt{H8}} (\Delta \tau_{\times}, \Delta \tau) \; \Xi \bigg]. 
    \label{ch3_discoTEX8_Wave_hyper}
\end{align} 
where, 
\begin{align}
    \textbf{TEXH8} = \bigg[ \textbf{I}  + \frac{\textbf{A}\cdot \textbf{A}}{42
    } \bigg],
    \label{discotex8_tex8}
\end{align}
and we have used the replacement, 
\begin{align}
    \textbf{U}^{(3)} &= \frac{d}{d\tau} \textbf{U}^{(2)} = \nonumber \\
    &= \textbf{L} \cdot \textbf{L} \cdot \textbf{L} \cdot \textbf{U} + \textbf{L} \cdot\textbf{L} \cdot \mathcal{S} + \textbf{L} \cdot \mathcal{S}^{(1)} + \mathcal{S}^{(2)}. 
    \label{discotex_h8_ddotuvector}
\end{align}
It will be further necessary to compute, 
\begin{align}
    g^{(3)}(\tau) &= \frac{d^{3}g(\tau)}{d\tau^{3}} = \frac{d}{d\tau}\big( g^{(2)}(\tau) \big) \nonumber \\
    &= \frac{1}{m!} J^{(3)}w^{m} - \frac{m}{m!}\big[3 J^{(2)} x^{(1)}_{p} + 3 J^{(1)} x^{(2)}_{p} + J x^{(3)}_{p} \big]w^{m-1} \nonumber \\
    &+ \frac{m(m-1)}{m!}\big[3 (x_{p}^{(1)})^{2} J^{(1)} + 3 x_{p}^{(1)} x_{p}^{(2)} \big]w^{m-2} \nonumber \\ 
    &- \frac{m(m-1)(m-2)}{m!} (x_{p}^{(1)})^{3} J w^{m-3}. 
    \label{discotex_h8_ddotg}
\end{align}

Similarly \texttt{DiscoTEX H10} is given by, 
\begin{align}
    &\textbf{U}_{n+1} = \textbf{U}_{n}+  \textbf{HFH10} \cdot  \bigg[ \textbf{A} \cdot  \bigg[ \textbf{TEXH10} \cdot \textbf{U}_{n} + \frac{ \Delta \tau}{9} \bigg(\textbf{s}_{n} - \textbf{s}_{n+1} \bigg) \nonumber \\
        &+\frac{\Delta \tau^{2}}{72} \bigg( \big(\textbf{s}^{(1)}_{n} +\textbf{s}^{(1)}_{n+1}\big) + \textbf{L} \cdot \big(\textbf{s}_{n} +\textbf{s}_{n+1}\big) \bigg) + \nonumber\\
        &+\frac{\Delta \tau^{3}}{1008} \bigg(\big(\textbf{s}^{(2)}_{n} -\textbf{s}^{(2)}_{n+1}\big) + \textbf{L} \cdot \big(\textbf{s}^{(1)}_{n} -\textbf{s}^{(1)}_{n+1}\big) + \textbf{L} \cdot \textbf{L} \cdot \big(\textbf{s}_{n} -\textbf{s}_{n+1}\big) \bigg) \nonumber \\ 
         &+\frac{\Delta \tau^{4}}{30240} \bigg(\big(\textbf{s}^{(3)}_{n} +\textbf{s}^{(3)}_{n+1}\big) + \textbf{L} \cdot \big(\textbf{s}^{(2)}_{n} +\textbf{s}^{(2)}_{n+1}\big) + \nonumber \\ 
         &+\textbf{L} \cdot \textbf{L} \cdot \big(\textbf{s}^{(1)}_{n} +\textbf{s}^{(1)}_{n+1}\big) + \textbf{L} \cdot \textbf{L} \cdot \textbf{L} \cdot \big(\textbf{s}_{n} +\textbf{s}_{n+1}\big) \bigg) \bigg]\nonumber \\
        &+  \frac{\Delta \tau}{2} \bigg( \textbf{s}_{n} + \textbf{s}_{n+1} \bigg) + \frac{\Delta \tau^{2}}{9} \bigg(\textbf{s}^{(1)}_{n} -\textbf{s}^{(1)}_{n+1}  \bigg) + \frac{\Delta \tau^{3}}{72} \bigg(\textbf{s}^{(2)}_{n} + \textbf{s}^{(2)}_{n+1}  \bigg)  \nonumber \\
        &+ \frac{\Delta \tau^{4}}{1008} \bigg(\textbf{s}^{(3)}_{n} -\textbf{s}^{(3)}_{n+1}  \bigg) + \frac{\Delta \tau^{5}}{30240} \bigg(\textbf{s}^{(4)}_{n} + \textbf{s}^{(4)}_{n+1}  \bigg) \nonumber \\
        &+ \Upsilon +  \textbf{J}_{\texttt{H10}} (\Delta \tau_{\times}, \Delta \tau) \; \Xi \bigg]. 
    \label{ch3_discoTEX8_Wave_hyper}
\end{align} 

where, 
\begin{align}
    \textbf{TEXH10} = \bigg[ \textbf{I}  + \textbf{A} \cdot \textbf{A} \cdot \bigg( \frac{\textbf{I}}{36} + \frac{\textbf{A}\cdot \textbf{A}}{1520}  \bigg)\bigg]
    \label{discotex10_tex10},
\end{align}
with the replacement 
\begin{align}
    \textbf{U}^{(4)} &= \frac{d}{d\tau} \textbf{U}^{(3)} = \nonumber \\
    &= \textbf{L} \cdot \textbf{L} \cdot \textbf{L} \cdot \textbf{L} \cdot \textbf{U} + \textbf{L} \cdot\textbf{L} \cdot \textbf{L} \cdot \mathcal{S} + \textbf{L} \cdot \textbf{L} \cdot \mathcal{S}^{(1)} \nonumber \\
    &+ \textbf{L} \cdot \mathcal{S}^{(2)} + \mathcal{S}^{(3)}. 
    \label{discotex_h10_ddotuvector}
\end{align}

In addition to equation \eqref{discotex_h8_ddotg}, the higher-order $g^{(4)}(\tau)$ will also be necessary, 
\begin{align}
    g^{(4)}(\tau) &= \frac{d^{4}g(\tau)}{d\tau^{4}} = \frac{d}{d\tau}\big( g^{(3)}(\tau) \big) \nonumber \\
     &=\frac{1}{m!}J^{(4)}(\tau) w^{m}(\tau) \nonumber \\
     &-\frac{m}{m!}\bigg[x^{(4)}_{p} J + 4 x^{(3)}_{p} J^{(1)} + 6 x^{(2)}_{p} J^{(2)}  + 4 x^{(1)}_{p} J^{(3)} \bigg]w^{m-1} \nonumber \\
     &+\frac{m(m-1)}{m!}\bigg[ 4 x^{(1)}_{p} x^{(3)}_{p} J + 12 x^{(1)}_{p} x^{(2)}_{p} J^{(1)} \nonumber \\ 
     &+ 6 (x^{(1)}_{p})^{2} J^{(1)} + 3 (x^{(2)}_{p})^{2} J \bigg]w^{m-2} \nonumber \\ 
     &- \frac{m(m-1)(m-2)}{m!} \bigg[ 6 (x^{(1)}_{p})^{2} x^{(2)}_{p} J +  4 (x^{(1)}_{p})^{3} J^{(1)} \bigg]w^{m-3} \nonumber\\
     &+ \frac{m(m-1)(m-2)(m-3)}{m!} (x_{p}^{(1)})^{4} J w^{m-4}. 
    \label{discotex_h10_dddotg}
\end{align}

Finally \texttt{DiscoTEX H12} is given by, 
\begin{align}
    &\textbf{U}_{n+1} = \textbf{U}_{n}+  \textbf{HFH12} \cdot  \bigg[ \textbf{A} \cdot  \bigg[ \textbf{TEXH12} \cdot \textbf{U}_{n} + \frac{ 5\Delta \tau}{44} \bigg(\textbf{s}_{n} + \textbf{s}_{n+1} \bigg) \nonumber \\
        &+\frac{\Delta \tau^{2}}{66} \bigg( \big(\textbf{s}^{(1)}_{n} +\textbf{s}^{(1)}_{n+1}\big) + \textbf{L} \cdot \big(\textbf{s}_{n} +\textbf{s}_{n+1}\big) \bigg) + \nonumber\\
        &+\frac{\Delta \tau^{3}}{792} \bigg(\big(\textbf{s}^{(2)}_{n} -\textbf{s}^{(2)}_{n+1}\big) + \textbf{L} \cdot \big(\textbf{s}^{(1)}_{n} -\textbf{s}^{(1)}_{n+1}\big) + \textbf{L} \cdot \textbf{L} \cdot \big(\textbf{s}_{n} -\textbf{s}_{n+1}\big) \bigg) \nonumber \\
      &+\frac{\Delta \tau^{4}}{15840} \bigg( \big(\textbf{s}^{(3)}_{n} +\textbf{s}^{(3)}_{n+1}\big) + \textbf{L} \cdot \big(\textbf{s}^{(2)}_{n} +\textbf{s}^{(2)}_{n+1}\big) + \textbf{L} \cdot \textbf{L} \cdot \big(\textbf{s}^{(1)}_{n} +\textbf{s}^{(1)}_{n+1}\big) \nonumber \\
      &+ \textbf{L} \cdot \textbf{L} \cdot \textbf{L} \cdot \big(\textbf{s}_{n} +\textbf{s}_{n+1}\big) \bigg) \nonumber \\
       &+\frac{\Delta \tau^{5}}{665280} \bigg( \big(\textbf{s}^{(4)}_{n} -\textbf{s}^{(4)}_{n+1}\big) + \textbf{L} \cdot \big(\textbf{s}^{(3)}_{n} -\textbf{s}^{(3)}_{n+1}\big) \nonumber \\ 
       &+ \textbf{L} \cdot \textbf{L} \cdot \big(\textbf{s}^{(2)}_{n} -\textbf{s}^{(2)}_{n+1}\big) + \textbf{L} \cdot \textbf{L} \cdot \textbf{L} \cdot \big(\textbf{s}^{(1)}_{n} -\textbf{s}^{(1)}_{n+1}\big) \nonumber \\
      &+ \textbf{L} \cdot \textbf{L} \cdot \textbf{L} \cdot \textbf{L} \cdot \big(\textbf{s}_{n} -\textbf{s}_{n+1}\big) \bigg) \bigg] \nonumber \\
    &+  \frac{\Delta \tau}{2} \bigg( \textbf{s}_{n} - \textbf{s}_{n+1} \bigg) + \frac{5\Delta \tau^{2}}{44} \bigg(\textbf{s}^{(1)}_{n} - \textbf{s}^{(1)}_{n+1}  \bigg) + \frac{\Delta \tau^{3}}{66} \bigg(\textbf{s}^{(2)}_{n} - \textbf{s}^{(2)}_{n+1}  \bigg) \nonumber \\
    &+\frac{\Delta \tau^{4}}{792} \bigg(\textbf{s}^{(3)}_{n} - \textbf{s}^{(3)}_{n+1}  \bigg) + \frac{\Delta \tau^{5}}{15840} \bigg(\textbf{s}^{(4)}_{n} - \textbf{s}^{(4)}_{n+1}  \bigg) \nonumber \\
    &+\frac{\Delta \tau^{6}}{665280} \bigg(\textbf{s}^{(5)}_{n} - \textbf{s}^{(5)}_{n+1}  \bigg) + \Upsilon +  \textbf{J}_{\texttt{H12}} (\Delta \tau_{\times}, \Delta \tau) \; \Xi \bigg]. 
    \label{ch3_discoTEX12_Wave_hyper}
\end{align} 
where, 
\begin{align}
    \textbf{TEXH12} = \bigg[ \textbf{I} + \textbf{A}\cdot\textbf{A}\cdot \bigg(\frac{\textbf{I}}{33} + \frac{\textbf{A}\cdot \textbf{A}}{7920} \bigg) \bigg],
    \label{discotex12_tex12}
\end{align}
with the replacement, 
\begin{align}
    \textbf{U}^{(5)} &= \frac{d}{d\tau} \textbf{U}^{(4)} = \nonumber \\
    &= \textbf{L} \cdot \textbf{L} \cdot \textbf{L} \cdot \textbf{L} \cdot \textbf{L} \cdot \textbf{U} + \textbf{L} \cdot \textbf{L} \cdot\textbf{L} \cdot \textbf{L} \cdot \mathcal{S} + \textbf{L} \cdot \textbf{L} \cdot  \textbf{L} \cdot \mathcal{S}^{(1)} \nonumber \\
    &+ \textbf{L} \cdot \textbf{L} \cdot \mathcal{S}^{(2)} +  \textbf{L} \cdot \mathcal{S}^{(3)} + \mathcal{S}^{(4)}, 
    \label{discotex_h12_ddotuvector}
\end{align}
and the additional higher-order $g^{(5)}(\tau)$ vector, 
\begin{align}
    g^{(5)}(\tau) &= \frac{d^{5}g(\tau)}{d\tau^{5}} = \frac{d}{d\tau}\big( g^{(4)}(\tau) \big) \nonumber \\
     &=\frac{1}{m!}J^{(5)}(\tau) w^{m}(\tau) \nonumber \\
     &-\frac{m}{m!}\bigg[J x^{(5)}_{p} + 5 x^{(4)}_{p} J^{(1)} + 10 x^{(3)} J^{(2)} + 10 x^{(2)} J^{(3)} + 5 x^{(1)} J^{(4)} \bigg]w^{m-1} \nonumber \\
     &+\frac{m(m-1)}{m!}\bigg[ 5 x^{(1)}_{p} x^{(4)}_{p} J + 10 x^{(2)}_{p} x^{(3)}_{p} J + 20 x^{(1)}_{p} x^{(3)}_{p} J^{(1)}   \nonumber \\
     &+15 (x^{(2)}_{p})^{2} J^{(1)} + 30 x^{(1)}_{p} x^{(2)}_{p} J^{(2)} + 10 (x^{(1)}_{p})^{2} J^{(3)}  \bigg]w^{m-2} \nonumber \\ 
     &- \frac{m(m-1)(m-2)}{m!} \bigg[10 (x_{p}^{(1)})^{2} x_{p}^{(3)} J + 15 x_{p}^{(1)} (x_{p}^{(2)})^{2} J \nonumber \\ 
     &+ 30 (x_{p}^{(1)})^{2}  x_{p}^{(2)} J^{(1)}  + 10 (x_{p}^{(1)})^{3} J^{(2)} \bigg]w^{m-3} \nonumber\\
     &+ \frac{m(m-1)(m-2)(m-3)}{m!} \bigg[10 (x^{(1)}_{p})^{3} x^{(2)}_{p} J + 5 (x^{(1)}_{p})^{4} J^{(1)}   \bigg]w^{m-4} \nonumber \\
     &- \frac{m(m-1)(m-2)(m-3)(m-4)}{m!} (x_{p}^{(1)})^{5} J w^{m-5}.
    \label{discotex_h12_dddotg}
\end{align}

\subsection{Time jumps for the implementation of \texttt{DiscoTEX} in the $(\tau,\sigma)$ \textit{minimal gauge} coordinate chart}\label{app_hyper_timejumps}
The time jumps necessary for the implementation $\textbf{J}_{\texttt{H2-12}}(\Delta t_{\times}, \Delta \tau )$ are explicitly given as: 
\begin{align}
    \label{hyper_wave_jb_0}
    \mathbb{J}_{0} &= \frac{4}{15} \cos{\tau_{c}} + \frac{272}{225}  \iu \sin{\tau_{c}}  \\
    \label{hyper_wave_Kb_0}
    \mathbb{K}_{0} &= \frac{4864}{3375} \iu \cos{\tau_{c}} -  \frac{128}{225} \sin{\tau_{c}} , \\
    \label{hyper_wave_Lb_0}
    \mathbb{L}_{0} &= - \frac{3136}{3375} \cos{\tau_{c}} - \frac{90368}{50625} \iu \sin{\tau_{c}}, \\
    \label{hyper_wave_Mb_0}
    \mathbb{M}_{0} &= - \frac{1724416 }{759375} \iu \cos{\tau_{c}} + \frac{69632}{50625} \sin{\tau_{c}}, \\
   \label{hyper_wave_Nb_0}  
    \mathbb{N}_{0} &= \frac{1475584}{759375} \cos{\tau_{c}} + \frac{33492992}{11390625} \iu \sin{\tau_{c}}, \\
    \label{hyper_wave_Ob_0}  
    \mathbb{O}_{0} &= \frac{657915904}{170859375} \iu \cos{\tau_{c}} - \frac{30507008}{11390625}\sin{\tau_{c}}, \\
    \label{hyper_wave_Pb_0}  
    \mathbb{P}_{0} &= - \frac{622084096}{170859375} \cos{\tau_{c}} - \frac{13014990848}{2562890625} \iu \sin{\tau_{c}}, \\
    \label{hyper_wave_Qb_0}  
    \mathbb{Q}_{0} &=  -\frac{258579890176}{38443359375} \iu \cos{\tau_{c}} + \frac{12585009152}{2562890625} \sin{\tau_{c}}, \\
    \label{hyper_wave_Rb_0}  
    \mathbb{R}_{0} &= \frac{253420109824}{38443359375} \cos{\tau_{c}} + \frac{5150958682112}{576650390625} \iu \sin{\tau_{c}}, \\
    \label{hyper_wave_Sb_0}  
    \mathbb{S}_{0} &=  \frac{102771504185344}{8649755859375} \iu \cos{\tau_{c}} - \frac{5089041317888}{576650390625} \sin{\tau_{c}}, \\
    \label{hyper_wave_Tb_0}  
    \mathbb{T}_{0} &= -\frac{102028495814656}{8649755859375} \cos{\tau_{c}} -\frac{2052458050224128 }{129746337890625} \iu \sin{\tau_{c}}, \\
    \label{hyper_wave_Ub_0}  
    \mathbb{U}_{0} &= -\frac{41013496602689536}{1946195068359375} \iu \cos{\tau_{c}} + \frac{2043541949775872}{129746337890625} \sin{\tau_{c}},
\end{align}

As we have done in Section of ref.\cite{da2024discotex} throughout equations ((145)-(159)) we need to incorporate the differential operators in $\textbf{L}$, as defined in Eqs.~(\eqref{ch34_red_1ode}, \eqref{ch34_hyperboloidal_L1_p} - \eqref{ch34_hyperboloidal_L2_p}) into the final time jump associated with the $\textbf{L}_{1,2}$ operators. We thus introduce the following additional higher-order terms,
\begin{align}
    \label{hyper_wave_jc_0}
    &\mathcal{J}_{0} =\mathcal{J}_{0,\chi} + \mathcal{J}_{0,\iota}  + \mathcal{J}_{0,\varepsilon} + \mathcal{J}_{0,\varrho},  \\
    \label{hyper_wave_Kc_0}
    &\mathcal{K}_{0} =\mathcal{K}_{0,\chi} + \mathcal{K}_{0,\iota}  + \mathcal{K}_{0,\varepsilon} + \mathcal{K}_{0,\varrho}, \\
    \label{hyper_wave_Lc_0}
    &\mathcal{L}_{0} = \mathcal{L}_{0,\chi} + \mathcal{L}_{0,\iota}  + \mathcal{L}_{0,\varepsilon} + \mathcal{L}_{0,\varrho}, \\
    \label{hyper_wave_Mc_0}
    &\mathcal{M}_{0} = \mathcal{M}_{0,\chi} + \mathcal{M}_{0,\iota}  + \mathcal{M}_{0,\varepsilon} + \mathcal{M}_{0,\varrho},\\
    &\mathcal{N}_{0} = \mathcal{N}_{0,\chi} + \mathcal{N}_{0,\iota}  + \mathcal{N}_{0,\varepsilon} + \mathcal{N}_{0,\varrho},\\
    &\mathcal{O}_{0} = \mathcal{O}_{0,\chi} + \mathcal{O}_{0,\iota}  + \mathcal{O}_{0,\varepsilon} + \mathcal{O}_{0,\varrho},\\
    &\mathcal{P}_{0} = \mathcal{P}_{0,\chi} + \mathcal{P}_{0,\iota}  + \mathcal{P}_{0,\varepsilon} + \mathcal{P}_{0,\varrho},\\
    &\mathcal{Q}_{0} = \mathcal{Q}_{0,\chi} + \mathcal{Q}_{0,\iota}  + \mathcal{Q}_{0,\varepsilon} + \mathcal{Q}_{0,\varrho},\\
    &\mathcal{R}_{0} = \mathcal{R}_{0,\chi} + \mathcal{R}_{0,\iota}  + \mathcal{R}_{0,\varepsilon} + \mathcal{R}_{0,\varrho},\\
    &\mathcal{S}_{0} = \mathcal{S}_{0,\chi} + \mathcal{S}_{0,\iota}  + \mathcal{S}_{0,\varepsilon} + \mathcal{S}_{0,\varrho},\\
    &\mathcal{T}_{0} = \mathcal{T}_{0,\chi} + \mathcal{T}_{0,\iota}  + \mathcal{T}_{0,\varepsilon} + \mathcal{T}_{0,\varrho},\\
    &\mathcal{U}_{0} = \mathcal{U}_{0,\chi} + \mathcal{U}_{0,\iota}  + \mathcal{U}_{0,\varepsilon} + \mathcal{U}_{0,\varrho},
\end{align}

Following ref.\cite{da2024discotex} we correct the differential operators which contain a coefficient given in terms of the $\sigma$ coordinate. In short\footnote{For more details see Appendices A and D of \cite{da2024discotex}.}, we define, 
\begin{align}
    &\mathcal{J}(\tau) =\mathcal{J}_{(0,\chi), \bar{m}} +
    \mathcal{J}_{(0,\iota),\bar{m}}  + \mathcal{J}_{(0,\varepsilon), \bar{m}} + \mathcal{J}_{(0,\varrho), \bar{m}}, 
\end{align}
where $\bar{m} = \{\bar{k}, \cdots, m - \bar{k}\}$ and $\bar{k}$ is variable depending on the coefficient of the differential operator in question.  Each of these terms is given as,
\begin{align}
    \label{wave_dx2_jtau_0_exp}
    &\mathcal{J}_{0, \chi} =  \sum^{\bar{m}}_{k=0} {\bar{m}\choose k} \tilde{\chi}^{(k)}(\xi_{p}) J_{\bar{m}-k}(\tau),\\
    \label{wave_dx_jtau_0_exp}
    &\mathcal{J}_{0, \iota} =  \sum^{\bar{m}}_{k=0} {\bar{m}\choose k} \tilde{\iota}^{(k)}(\xi_{p}) J_{\bar{m}-k}(\tau),\\
    \label{wave_dx_jtau_0_exp}
   &\mathcal{J}_{0, \varepsilon} =  \sum^{\bar{m}}_{k=0} {\bar{m}\choose k} \tilde{\varepsilon}^{(k)}(\xi_{p}) \dot{J}_{\bar{m}-k}(\tau),\\
    \label{wave_jtau_0_exp}
    &\mathcal{J}_{0, \varrho} =  \sum^{\bar{m}}_{k=0} {\bar{m}\choose k} \tilde{\varrho}^{(k)}(\xi_{p}) \dot{J}_{\bar{m}-k}(\tau), \\
    \label{jjtau0}
    &\mathcal{J}_{0}(\tau) =  \mathcal{J}_{0, \chi}|_{\bar{k}=2} + \mathcal{J}_{0, \iota}|_{\bar{k}=1} + \mathcal{J}_{0, \varepsilon}|_{\bar{k}=0} + \mathcal{J}_{0, \varrho}|_{\bar{k}=1}.
\end{align}
All the remainder higher-order terms are determined in an anologous trivial process which we omit for brevity. Below we give the final results:
\begin{align}
    \label{wave_jc_0_exp}
    \mathcal{J}_{0} &= \frac{4\big( \big(225 + 1216\iu \big) \cos{\tau_{c}} - \big(480-1020\iu \big) \sin{\tau_{c}} \big)}{3375}, \\
    \label{wave_Kc_0_exp}
    \mathcal{K}_{0} &= \frac{64\big( \big(-735 + 1140\iu \big) \cos{\tau_{c}} - \big(450-1412\iu \big) \sin{\tau_{c}} \big)}{50625},\\
    \label{wave_Lc_0_exp}
    \mathcal{L}_{0} &= -\frac{64\big(\big(11025 + 26944\iu \big) \cos{\tau_{c}} - \big(16320-21180\iu \big) \sin{\tau_{c}} \big)}{759375},\\
    \label{wave_Mc_0_exp}
    \mathcal{M}_{0} &= \frac{(22133760 - 25866240\iu)\cos{\tau_{c}}}{11390625} \nonumber \\
    &+ \frac{(15667200 + 33492992 \iu) \sin{\tau_{c}}}{11390625}, \\
    \label{wave_Nc_0_exp}
    \mathcal{N}_{0} &= \frac{1024(324225 + 642496\iu)\cos{\tau_{c}}}{170859375} \nonumber \\
    &-\frac{1024 (446880 - 490620\iu) \sin{\tau_{c}}}{170859375},\\
    \label{wave_Oc_0_exp}
    \mathcal{O}_{0} &= \frac{16384(-569535 + 602340\iu)\cos{\tau_{c}}}{2562890625} \nonumber \\
    &-\frac{16384(418950 + 794372\iu)\sin{\tau_{c}}}{2562890625}, \\
    \label{wave_Pc_0_exp}
    \mathcal{P}_{0} &= \frac{16384(8543025 + 15782464\iu)\cos{\tau_{c}}}{38443359375} \nonumber \\
    &- \frac{16384(11521920 - 11915580\iu) \sin{\tau_{c}}}{38443359375}\\
    \label{wave_Qc_0_exp}
    \mathcal{Q}_{0} &= \frac{(3801301647360 - 3878698352640\iu)\cos{\tau_{c}}}{576650390625} \nonumber \\
    &+ \frac{ (2831627059200 + 5150958682112\iu) \sin{\tau_{c}}}{576650390625} \\
    \label{wave_Rc_0_exp}
    \mathcal{R}_{0} &= \frac{262144(217512225 + 392042176\iu)\cos{\tau_{c}}}{8649755859375} \nonumber \\
    &- \frac{262144(291197280 - 294740220\iu) \sin{\tau_{c}}}{8649755859375} \\
    \label{wave_Sc_0_exp}
    \mathcal{S}_{0} &= \frac{4194304(-364882335 + 367539540\iu)\cos{\tau_{c}}}{129746337890625} \nonumber \\
    &- \frac{4194304(272997450 + 489344132\iu) \sin{\tau_{c}}}{129746337890625} \\
    \label{wave_Tc_0_exp}
    \mathcal{T}_{0} &= -\frac{4194304(5473235025 + 9778379584\iu)\cos{\tau_{c}}}{1946195068359375} \nonumber \\
    &-\frac{4194304(7308275520-7340161980\iu) \sin{\tau_{c}}}{1946195068359375} \\
    \label{wave_Vc_0_exp}
    \mathcal{U}_{0} &= \frac{(613597550959656960 - 615202449040343040\iu)\cos{\tau_{c}}}{29192926025390625}\nonumber \\
    &+ \frac{(459796938699571200 + 819841959232274432\iu) \sin{\tau_{c}}}{29192926025390625} 
\end{align}

Finally, we note all these equations and others are given before taking the limit as previously explained.


\end{document}